\documentclass[12pt,reqno]{article}

\usepackage[usenames]{color}
\usepackage[dvipsnames]{xcolor}

\usepackage{graphicx}
\usepackage{amscd}

\definecolor{webgreen}{rgb}{0,.5,0}
\definecolor{webbrown}{rgb}{.6,0,0}

\usepackage{color}
\usepackage{fullpage}
\usepackage{float}

\usepackage{graphics,amsmath,amssymb}
\usepackage{amsthm}
\usepackage{amsfonts}
\usepackage{latexsym}
\usepackage{epsf}
\usepackage{mathrsfs}

\usepackage[colorlinks=true,
linkcolor=webgreen,
filecolor=webbrown,
citecolor=webgreen]{hyperref} 

\def\modd#1 #2{#1\ \mbox{\rm (mod}\ #2\mbox{\rm )}}

\newcommand{\seqnum}[1]{\href{http://oeis.org/#1}{\underline{#1}}}

\begin{document}

\theoremstyle{plain}
\newtheorem{theorem}{Theorem}
\newtheorem{corollary}[theorem]{Corollary}
\newtheorem{lemma}[theorem]{Lemma}
\newtheorem{proposition}[theorem]{Proposition}

\theoremstyle{definition}
\newtheorem{definition}[theorem]{Definition}
\newtheorem{example}[theorem]{Example}
\newtheorem{conjecture}[theorem]{Conjecture}

\theoremstyle{remark}
\newtheorem{remark}[theorem]{Remark}

\begin{center}
\vskip 1cm{\LARGE\bf  Areas Between Cosines}

\vskip 1cm
\large

Muhammad Adam Dombrowski\\
Pennsbury High School\\
Fairless Hills, PA 19030\\
 USA\\
\href{mailto:muhammadadamdombrowski@gmail.com}{\tt muhammadadamdombrowski@gmail.com}\\ 

\ \\

Gregory Dresden\\ 
Washington \& Lee University \\
Lexington, VA 24450 \\
USA\\
\href{mailto:dresdeng@wlu.edu}{\tt dresdeng@wlu.edu} \\
\end{center}
\vskip0.2in

\begin{abstract} 
We find the area between $\cos^n x$ and $\cos^n kx$ over the interval $[0,\pi]$ as $k$ heads to infinity. We establish recursive formulas for these areas, and we show that these  areas are related to the   coefficients from two exponential generating functions involving   $\arcsin x$.
\end{abstract}

\section{Introduction}

\begin{figure}[h]
\begin{center}	\includegraphics[width=3.5in]{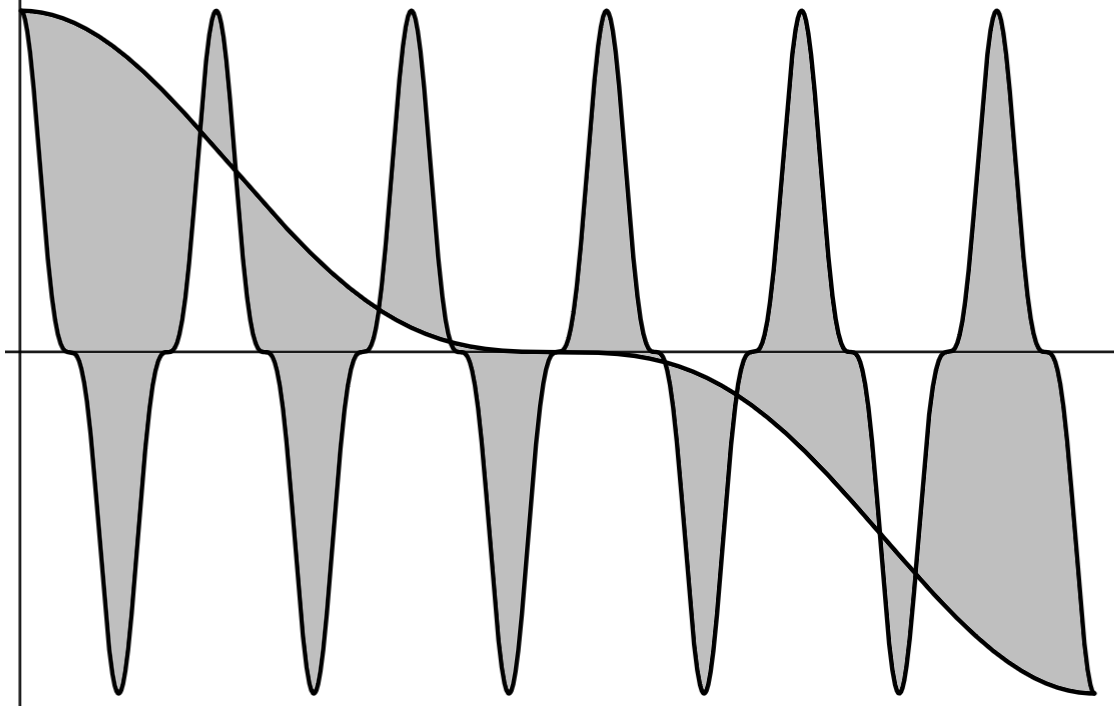}
\end{center}
\caption{Region between $\cos^3 x$ and $\cos^3 11x$.}\label{f.firstpic}
\end{figure}

Figure \ref{f.firstpic} shows the region between 
$\cos^3 x$ and $\cos^3 11x$ over 
the interval $[0,\pi]$.

It is not hard to directly calculate the area of this region. We could write it as
\[
\frac{6}{33}\left( \cot \frac{3\pi}{12} + 
9\cot \frac{\pi}{12}\right) -
\frac{5}{33}\left( \cot \frac{3\pi}{10} + 
9\cot \frac{\pi}{10}\right) \approx 1.981887\dots, 
\]
which has a pleasant symmetry, or we could write it as 
\[
\frac{1}{33}\left( 144 + 54 \sqrt{3} - 
5\left(19+9\sqrt{5}\right) \sqrt{5-2\sqrt{5}}\right) \approx 1.981887\dots,
\]
which has a purely algebraic expression. 
These are both interesting answers, but in this article we  are more interested in what happens when we replace the 
11 in $\cos^3 11 x$ with
$k$ (giving us $\cos^3 kx$) and then 
find the area between $\cos^3 x$ and $\cos^3 kx$ as  $k$ goes to infinity. To be precise, we want to find 
\[
\lim_{k \to \infty} \int_0^\pi \Big| \cos^3 x - \cos^3 kx \Big| \, dx.
\]
In this case, the limiting area turns out to be 
\[
\lim_{k \to \infty} \int_0^\pi \Big| \cos^3 x - \cos^3 kx \Big| \, dx = \frac{56}{9 \pi} \approx 1.980594\dots
\] 
which is rather surprising (and also fairly close numerically to the two expressions above). 

Of course, there is no reason to restrict ourselves to just looking at the {\em third} power of cosine. With this in mind, we define $A_n$ to be the limiting area (as $k \to \infty$) of the region between $\cos^n x$ and $\cos^n kx$ over 
the interval $[0,\pi]$.
In other words, we define
\begin{equation}\label{e.An.definition}
A_n = \lim_{k \to \infty} \int_0^\pi \Big| \cos^n x - \cos^n kx \Big| \, dx.
\end{equation}

In what follows, we will find formulas for $A_n$ involving sums with binomial coefficients (Theorems \ref{t.oddpower} and \ref{t.evenpower}).
We then find recursive formulas for $A_n$ involving just $A_{n-2}$ (Theorem \ref{t.recursion}). Finally,  in Theorem 
\ref{t.conjecture}
we establish formulas for $A_n$ that involve sums with double factorials, and we connect these formulas with two entries
in the 
On-Line Encyclopedia of Integer Sequences (OEIS) \cite{oeis}. These two entries are 
\seqnum{A296726}  and 
\seqnum{A372324},  and they are 
related to the exponential 
generating functions for $\arcsin x/(1-x)$ and $\arcsin^2 x/(2(1-x))$, respectively.

\section{Area formulas} 
Since $\cos^n x$ looks quite different on the interval $[0,\pi]$ depending on the parity of $n$ (as seen by comparing Figures \ref{f.311color} and \ref{f.47color}), it is reasonable to separate our discussion of areas into the following two cases. 
\subsection{Odd $n$}
To set the stage, we begin with a few values for the limiting area $A_n$  when $n$ is odd:
\begin{align*}
    A_1 &=   \frac{8}{\pi} = \displaystyle\frac{8\cdot 1}{(1)^2 \pi}, 
        &  
    A_3 &=  \frac{56}{9 \pi} = \displaystyle\frac{8\cdot 7}{(1\cdot 3)^2\pi}, 
        \\[2.5ex]
    A_5 &= \frac{1192}{225\pi} =  \displaystyle\frac{8\cdot 149}{(1\cdot 3\cdot 5)^2\pi},
        &
 A_7 &= \frac{17228}{3675 \pi} = \displaystyle\frac{8\cdot 6483}{(1\cdot 3\cdot 5\cdot 7)^2\pi}.
\end{align*}
What can we say about these numbers
$1, 7, 149, 6483$ that appear in the numerators of $A_n$ for $n$ odd? We will show that these are equal to every other term in the sequence
\seqnum{A296726} 
in the OEIS, where we learn that they also appear as  coefficients in the exponential generating function for 
$\arcsin x/(1-x)$. See Theorem \ref{t.conjecture} for details. 

We also note that the value $A_1 = 8/\pi$ was
the answer to Problem 2191 which appeared in a recent issue 
\cite{2191} of  {\em Mathematics Magazine}. 

Here is our first theorem.

\begin{theorem}\label{t.oddpower}
With  $A_n$ as defined above in equation 
(\ref{e.An.definition}), then 
\begin{equation}\label{e.firstAp}
\text{for $n$ odd,} \qquad A_n = 
\frac{8}{2^{n-1}\pi}  \sum_{j=0}^{(n-1)/2} \binom{n}{j}\frac{1}{(n-2j)^2}.
\end{equation}
\end{theorem}

To prove Theorem \ref{t.oddpower}, we will need several technical results that we present and prove in Section \ref{s.technical}. The proof of Theorem \ref{t.oddpower} is in 
Section \ref{s.proof1}

\subsection{Even $n$}
Next, we present a few values for the limiting area $A_n$ for when $n$ is even:
\begin{align*}
    A_2 &= \frac{4}{\pi} =  \frac{16\cdot 1}{(2)^2\pi}, 
        &  
    A_4 &= \frac{4}{\pi} = \frac{16\cdot 16}{(2\cdot 4)^2\pi}, 
        \\[2.5ex]
    A_6 &= \frac{34}{9\pi} = \frac{16\cdot 544}{(2\cdot 4\cdot 6)^2\pi},
        &
 A_8 &= \frac{32}{9 \pi} = \frac{16\cdot 32768}{(2\cdot 4\cdot 6\cdot 8)^2\pi}.
\end{align*}

What can we say about these numbers
$1, 16, 544, 32768, \dots$ that appear in the numerators of $A_n$ for $n$ even? 
We will show that these are equal to 
alternate terms in the sequence
\seqnum{A372324} 
in the OEIS, where we learn that they also appear as  coefficients in the exponential generating function for 
$\arcsin^2 x/(2(1-x))$. See Theorem \ref{t.conjecture} for details. 

For now, we have the following
theorem.

\begin{theorem}\label{t.evenpower}
With  $A_n$ as defined above in equation (\ref{e.An.definition}), then 
\begin{align}
\text{for $n \equiv \modd{2} {4}$,} \qquad A_n &= \frac{16}{2^{n}\pi} \sum_{j=0}^{(n-2)/4} \binom{n}{2j}\frac{1}{(n/2-2j)^2}, \label{e.Apforp24} \\
\intertext{and }
\text{for $n \equiv \modd{0} {4}$,} \qquad A_n &=  \frac{16}{2^{n}\pi} \sum_{j=0}^{(n-4)/4} \binom{n}{2j+1}\frac{1}{(n/2-(2j+1))^2}. \label{e.Apforp04}
\end{align}

\end{theorem}

Just as with Theorem \ref{t.oddpower},  the proof of 
Theorem \ref{t.evenpower} will require some technical results that we present and prove in Section \ref{s.technical}. The proof of Theorem \ref{t.evenpower} is in 
Section \ref{s.proof2}

\section{Recursion formulas}
Thanks to Theorems \ref{t.oddpower} and \ref{t.evenpower}, we can produce the following rather simple recursion formula.
\begin{theorem}\label{t.recursion}
With  $A_n$ as defined above, then with $n\geq 3$ we have
\begin{displaymath}
A_n = \frac{n-1}{n}A_{n-2} + 
\begin{cases}
\displaystyle \frac{8}{n^2\pi}, &  \text{for $n$ odd;}
   \\[2.75ex]
 \displaystyle\frac{16}{n^2\pi}, &  
            \text{for $n$ even.}
\end{cases}
\end{displaymath}
\end{theorem}

\begin{proof}
We start with the easily-verified statement that
\[
\binom{n}{j} \cdot (n-j)\cdot j = 
\binom{n-2}{j-1} \cdot (n-1)\cdot n.
\]
Next, we multiply both sides by 4 to get 
\[
\binom{n}{j} \cdot (2n-2j)\cdot 2j = 
4 \cdot \binom{n-2}{j-1} \cdot (n-1)\cdot n.
\]
We now divide both sides by $(n-2j)^2= ((n-2) - 2(j-1))^2$, to get
\[
\binom{n}{j} \cdot \frac{(2n-2j)\cdot 2j}{(n-2j)^2} = 
4 \cdot \binom{n-2}{j-1} \cdot 
\frac{(n-1)\cdot n}{((n-2)-2(j-1))^2}.
\]
In the numerator on the left, we first replace
$2n-2j$ with $n+n-2j$ and then  replace 
$2j$ with $n-(n-2j)$, giving us 

\[
\binom{n}{j} \cdot \frac{(n+(n-2j))\cdot (n-(n-2j))}{(n-2j)^2} = 
4 \cdot \binom{n-2}{j-1} \cdot 
\frac{(n-1)\cdot n}{((n-2)-2(j-1))^2}.
\]
We multiply out the numerator on the left to get 
\[
\binom{n}{j} \cdot \frac{n^2 - (n-2j)^2}{(n-2j)^2} = 
4 \cdot \binom{n-2}{j-1} \cdot 
\frac{(n-1)\cdot n}{((n-2)-2(j-1))^2},
\]
and a further simplification on the left gives us 
\[
\binom{n}{j} \cdot \frac{n^2}{(n-2j)^2} 
- \binom{n}{j} = 
4 \cdot \binom{n-2}{j-1} \cdot 
\frac{(n-1)\cdot n}{((n-2)-2(j-1))^2}.
\]
We now divide both sides by $n^2$ and simplify a bit more to obtain
\begin{equation}\label{e.mainbinomial}
\binom{n}{j}  \frac{1}{(n-2j)^2} 
- \frac{1}{n^2}\binom{n}{j} = 
 \frac{4(n-1)}{n} \binom{n-2}{j-1}
\frac{1}{((n-2)-2(j-1))^2}.
\end{equation}

At this point, we will consider the three cases
of $n$ odd, $n \equiv \modd{2} {4}$, and $n \equiv \modd{0} {4}$. 

{\bf  First, we suppose $n$ is odd}. We sum both
sides of equation (\ref{e.mainbinomial}) from 
$j=1$ to $j=(n-1)/2$ to obtain
\[
\sum_{j=1}^{(n-1)/2}\binom{n}{j}  \frac{1}{(n-2j)^2} 
- \frac{1}{n^2}\binom{n}{j} = 
  \frac{4(n-1)}{n}
 \sum_{j=1}^{(n-1)/2}\binom{n-2}{j-1}  
\frac{1}{((n-2)-2(j-1))^2}.
\]
On the left, we can start that sum at $j=0$ instead of $j=1$ without changing the value. On the right, we re-index the sum by using $j'= j-1$, so that $j'$ starts at $j'=0$ and ends at $j' = (n-3)/2$. After distributing the sum on the left, this gives us
\[
\sum_{j=0}^{(n-1)/2}\binom{n}{j}  \frac{1}{(n-2j)^2} 
- \frac{1}{n^2}\sum_{j=0}^{(n-1)/2}\binom{n}{j} = 
  \frac{4(n-1)}{n}
 \sum_{j'=0}^{(n-3)/2}\binom{n-2}{j'}  
\frac{1}{((n-2)-2j')^2}.
\]
Thanks to our value for $A_n$ in equation
(\ref{e.firstAp}) for $n$ odd, we can re-write the 
above equation as 
\[
 A_n \frac{2^{n-1} \pi}{8}
- \frac{1}{n^2}\sum_{j=0}^{(n-1)/2}\binom{n}{j} = 
  \frac{4(n-1)}{n}\cdot A_{n-2}\frac{2^{(n-2)-1} \pi}{8}. 
\]
Since $n$ is odd, the sum on the left of the above equation is exactly half of the entire sum of the $n$th row of Pascal's triangle. The entire sum would be $2^n$, so the sum in the above equation would be $2^{n-1}$. 
So, after combining the 4 and the $2^{(n-2)-1}$ on the right, we have
\[
 A_n \frac{2^{n-1} \pi}{8}
- \frac{1}{n^2}2^{n-1} = 
  \frac{(n-1)}{n}\cdot A_{n-2}\frac{2^{n-1} \pi}{8}. 
\]
If we now divide everything
by ${2^{n-1} \pi}/{8}$ and re-arrange the terms, we obtain
\[
A_n = \frac{n-1}{n} A_{n-2} + \frac{8}{\pi n^2},
\]
as desired (for $n$ odd). 

{\bf Next, we consider} $n \equiv \modd{2} {4}$. Looking back at equation (\ref{e.mainbinomial}), we 
replace $j$ with $2j$, giving us 
\begin{equation*}
\binom{n}{2j}  \frac{1}{(n-4j)^2} 
- \frac{1}{n^2}\binom{n}{2j} = 
 \frac{4(n-1)}{n} \binom{n-2}{2j-1}
\frac{1}{((n-2)-2(2j-1))^2}.
\end{equation*}
We factor out $2^2$ from the $(n-4j)^2$ in the denominator on the left, and likewise from the denominator on the right, giving us 
\begin{equation*}\label{e.mainbinomial2j.2}
\binom{n}{2j}  \frac{1}{4(n/2-2j)^2} 
- \frac{1}{n^2}\binom{n}{2j} = 
 \frac{4(n-1)}{n} \binom{n-2}{2j-1}
\frac{1}{4((n-2)/2-(2j-1))^2}.
\end{equation*}
We multiply through by 4 to get 
\begin{equation}\label{e.mainbinomial2j.3}
\binom{n}{2j}  \frac{1}{(n/2-2j)^2} 
- \frac{4}{n^2}\binom{n}{2j} = 
 \frac{4(n-1)}{n} \binom{n-2}{2j-1}
\frac{1}{((n-2)/2-(2j-1))^2}.
\end{equation}
We sum both
sides of equation (\ref{e.mainbinomial2j.3}) from 
$j=1$ to $j=(n-2)/4$ to obtain
\[
\sum_{j=1}^{(n-2)/4}\binom{n}{2j}  \frac{1}{(n/2-2j)^2} 
- \frac{4}{n^2}\binom{n}{2j} = 
  \frac{4(n-1)}{n}
 \sum_{j=1}^{(n-2)/4}\binom{n-2}{2j-1}  
\frac{1}{((n-2)/2-(2j-1))^2}.
\]
On the left, we can start that sum at $j=0$ instead of $j=1$ without changing the value. On the right, we re-index the sum by using $j'= j-1$, so that $j'$ starts at $j'=0$ and ends at $j' = (n-6)/4$. This gives us
\[
\sum_{j=0}^{(n-2)/4}\binom{n}{2j}  \frac{1}{(n/2-2j)^2} 
- \frac{4}{n^2}\binom{n}{2j} = 
  \frac{4(n-1)}{n}
 \sum_{j'=0}^{(n-6)/4}\binom{n-2}{2j'+1}  
\frac{1}{((n-2)/2-(2j'+1))^2}.
\]
Thanks to our value for $A_n$ in equation
(\ref{e.Apforp24}) for $n\equiv \modd{2} {4}$, we can
re-write the first sum in the above equation as 
$\displaystyle A_n \cdot 2^n \pi/16$. When we do so, it gives us 
\[
A_n \frac{2^n \pi}{16}
- \frac{4}{n^2}\sum_{j=0}^{(n-2)/4}\binom{n}{2j} = 
  \frac{4(n-1)}{n}
 \sum_{j'=0}^{(n-6)/4}\binom{n-2}{2j'+1}  
\frac{1}{((n-2)/2-(2j'+1))^2}.
\]
Likewise, since $n \equiv \modd{2} {4}$, then 
$n-2 \equiv \modd{0} {4}$, and so if we use equation 
(\ref{e.Apforp04}) for $A_{n-2}$ then we 
recognize that the sum on the right of the above equation is equal to 
$\displaystyle A_{n-2}\cdot (2^{n-2}\pi)/16$.
This means we can now re-write the above equation as \[
A_n \frac{2^n \pi}{16}
- \frac{4}{n^2}\sum_{j=0}^{(n-2)/4}\binom{n}{2j} = 
  \frac{4(n-1)}{n}
 A_{n-2}\frac{2^{n-2}\pi}{16}.
\]
Since $n \equiv \modd{2} {4}$, then the sum on the left of the above equation is exactly one quarter of the entire sum of the $n$th row of Pascal's triangle. The entire sum would be $2^n$, so the sum in the above equation would be 
$2^{n-2}$, giving us
\[
A_n \frac{2^n \pi}{16}
- \frac{4}{n^2}2^{n-2} = 
  \frac{4(n-1)}{n}
 A_{n-2}\frac{2^{n-2}\pi}{16}.
\]
After combining the $4$ and the $2^{n-2}$ on the left and on the right, we have 
\[
A_n \frac{2^n \pi}{16}
- \frac{1}{n^2}2^n = 
  \frac{(n-1)}{n}
 A_{n-2}\frac{2^n\pi}{16}.
\]
If we now divide everything by 
${2^n \pi}/{16}$ and re-arrange the terms, we obtain
\[
A_n = \frac{n-1}{n} A_{n-2} + \frac{16}{n^2\pi},
\]
as desired (for $n \equiv \modd{2} {4}$).

{\bf Finally, we consider} $n \equiv \modd{0} {4}$.
 Looking back once more at equation (\ref{e.mainbinomial}), we  
 first factor out $2^2$ from the $(n-2j)^2$ 
 in the denominator on the left, and likewise from the denominator in the right. We also
 replace $n-2$ with $q$ in the expression on the right, leaving us with 
\[
\binom{n}{j}  \frac{1}{4(n/2-j)^2} 
- \frac{1}{n^2}\binom{n}{j} = 
 \frac{4(n-1)}{n} \binom{q}{j-1}
\frac{1}{4(q/2-(j-1))^2}.
\]
We now multiply through by 4, and 
replace $j$ with $2j+1$, giving us 
\begin{equation}\label{e.mainbinomial2j0}
\binom{n}{2j+1}  \frac{1}{(n/2-(2j+1))^2} 
- \frac{4}{n^2}\binom{n}{2j+1} = 
 \frac{4(n-1)}{n} \binom{q}{2j}
\frac{1}{(q/2-2j)^2}.
\end{equation}
We sum both
sides of equation (\ref{e.mainbinomial2j0}) from 
$j=0$ to $j=(n-4)/4$ to obtain
\[
\sum_{j=0}^{(n-4)/4}\binom{n}{2j+1}  \frac{1}{(n/2-(2j+1))^2} 
- \frac{4}{n^2}\binom{n}{2j+1} = 
  \frac{4(n-1)}{n}
 \sum_{j=0}^{(n-4)/4}\binom{q}{2j}  
\frac{1}{(q/2-2j)^2}.
\]
Thanks to our value for $A_n$ in equation
(\ref{e.Apforp04}) for $n\equiv \modd{0} {4}$, we can re-write the  sum of the first expression on the left above as $\displaystyle A_n\cdot(2^n \pi)/16$. When we do so (after 
replacing $n-4$ with $q-2$ in the upper bound of the sum on the right) it give us 
\[
A_n \frac{2^n \pi}{16}
- \frac{4}{n^2}\sum_{j=0}^{(n-4)/4}\binom{n}{2j+1} = 
  \frac{4(n-1)}{n}
 \sum_{j=0}^{(q-2)/4}\binom{q}{2j}  
\frac{1}{(q/2-2j)^2}.
\]
Since $n \equiv \modd{0} {4}$ and since $q=n-2$,
then 
$q \equiv \modd{2} {4}$, and so if we use equation 
(\ref{e.Apforp24}) for $A_{q}=A_{n-2}$ then we 
recognize that the sum on the right is equal to 
$\displaystyle A_{n-2}\cdot(2^{n-2} \pi)/16$.
This means we can re-write the above equation as 
\[
A_n \frac{2^n \pi}{16}
- \frac{4}{n^2}\sum_{j=0}^{(n-4)/4}\binom{n}{2j+1} = 
  \frac{4(n-1)}{n}
 A_{n-2}\frac{2^{n-2} \pi}{16}.
\]
Since $n \equiv \modd{0} {4}$, then the sum on the left of the above equation is exactly one quarter of the entire sum of the $n$th row of Pascal's triangle. The entire sum would be $2^n$, so the sum in the above equation would be $2^{n-2}$, giving us
\[
A_n \frac{2^n \pi}{16}
- \frac{4}{n^2}2^{n-2} = 
  \frac{4(n-1)}{n}
 A_{n-2}\frac{2^{n-2} \pi}{16}.
\]
After combining the $4$ and the $2^{n-2}$ on the left and on the right, we have 
\[
A_n \frac{2^n \pi}{16}
- \frac{1}{n^2}2^n = 
  \frac{(n-1)}{n}
 A_{n-2}\frac{2^n\pi}{16}.
\]
If we now divide everything by 
${2^n \pi}/{16}$ and re-arrange the terms, we obtain
\[
A_n = \frac{n-1}{n} A_{n-2} + \frac{16}{n^2\pi},
\]
as desired (for $n \equiv \modd{0} {4}$). 

Having covered all the cases for $n$, this completes the proof. 
\end{proof}

\section{Connections to the OEIS}

As we mentioned earlier, the numbers that appear in Theorem \ref{t.oddpower} 
are related to the sequence \seqnum{A296726},
and likewise those in Theorem \ref{t.evenpower}
appear in sequence \seqnum{A372324}. 
Here is the connection.

\begin{theorem}\label{t.conjecture}
For $A_n$ as defined above in equation (\ref{e.An.definition}), then 
\begin{align}
    \text{for $n$ odd,} \qquad A_n &= \frac{8}{\pi}\cdot \frac{n!}{(n!!)^2} \cdot \sum_{{j}=0}^{(n-1)/2} \frac{(2{j}-1)!!}{(2{j})!!} \frac{1}{2{j}+1}, \label{e.secondAp} \\[2.0ex]
\text{and for $n$ even,} \qquad A_n &= \frac{16}{\pi}\cdot \frac{n!}{(n!!)^2} \cdot \sum_{{j}=0}^{(n-2)/2} \frac{(2{j})!!}{(2{j}+1)!!} \frac{1}{2{j}+2}. \label{e.secondApeven}
\end{align}
Furthermore, the numbers 
\begin{equation}\label{e.A296726}
{n!}  \sum_{{j}=0}^{(n-1)/2} \frac{(2{j}-1)!!}{(2{j})!!} \frac{1}{2{j}+1} \qquad \text{for $n$ odd,}
\end{equation}
from equation (\ref{e.secondAp}) above, appear as every other entry in 
{\rm \seqnum{A296726}}, the terms from the exponential generating function for 
$\arcsin x/(1-x)$. Likewise, 
the numbers 
\begin{equation}\label{e.Anotyet}
{n!}  \sum_{{j}=0}^{(n-2)/2} \frac{(2{j})!!}{(2{j}+1)!!} \frac{1}{2{j}+2} \qquad \text{for $n$ even,}
\end{equation}
from equation (\ref{e.secondApeven}) above, appear as every other entry in 
{\rm \seqnum{A372324}}, the terms from 
exponential generating function for 
$\arcsin^2 x/(2(1-x))$.
\end{theorem}

We recall that  the notation $n!$ refers to the usual factorial function, and the notation 
$n!!$ is the less-familiar {\em double factorial} function \cite{Gould} although, to be honest, $n!!$ should be called an ``every other factorial"  instead. Here is the definition:
\begin{align*}
(2{j})!! &= (2{j})(2{j}-2)(2{j}-4) \cdots 6\cdot 4 \cdot 2,\\
(2{j}+1)!! &= (2{j}+1)(2{j}-1)(2{j}-3) \cdots 5\cdot 3 \cdot 1.
\end{align*}
We also agree that $0!! = (-1)!! = 1$. The sequence of double factorials is 
\seqnum{A006882}.

\begin{proof}[Proof of Theorem \ref{t.conjecture}] 
{\bf We begin with $n$ odd}. We   define $A_n'$ to be the right-hand side of equation (\ref{e.secondAp}), so  that 
\begin{equation}\label{e.bothsidesApprime}
 A_n' =  \frac{8}{\pi}\cdot \frac{n!}{(n!!)^2} \cdot \sum_{{j}=0}^{(n-1)/2} \frac{(2{j}-1)!!}{(2{j})!!} \frac{1}{2{j}+1}.
\end{equation}
We know from  Theorem \ref{t.recursion} that for $n$ odd, we have
\[
A_n = \frac{n-1}{n}A_{n-2} + 
\frac{8}{n^2\pi}.
\]
We now seek to
prove that 
\begin{equation}\label{e.ApApp}
A_n' = \frac{n-1}{n}A_{n-2}' + 
\frac{8}{n^2\pi}
\end{equation}
This, along with the fact that $A_1 = A_1' = 8/\pi$, is all we will need to conclude that $A_n = A_n'$ for all
$n$, thus proving the validity of equation 
(\ref{e.secondAp}). 

Now, going back to equation  (\ref{e.bothsidesApprime}), we replace $n$ with $n-2$ to give us 
\[
A_{n-2}' 
= 
 \frac{8}{\pi} \cdot \frac{(n-2)!}{((n-2)!!)^2} \cdot \sum_{{j}=0}^{(n-3)/2}\frac{(2{j}-1)!!}{(2{j})!!}\frac{1}{2{j}+1}.
\]
Next, we multiply both sides by $(n-1)/n$ and then add $8/(n^2 \pi)$ to give us 
\[
\frac{n-1}{n}A_{n-2}' + \frac{8}{n^2\pi}  
= 
\frac{8}{\pi} \left( \frac{1}{n^2} + \frac{n-1}{n} \cdot \frac{(n-2)!}{((n-2)!!)^2} \sum_{{j}=0}^{(n-3)/2}\frac{(2{j}-1)!!}{(2{j})!!}\frac{1}{2{j}+1}  \right).
\]
Since 
\begin{equation}\label{e.easybefore}
\frac{n-1}{n} \cdot \frac{(n-2)!}{((n-2)!!)^2} = 
\frac{n(n-1)}{n^2} \cdot \frac{(n-2)!}{((n-2)!!)^2}
= 
\frac{n!}{(n!!)^2},
\end{equation}
we can re-write the previous equation as 
\[
\frac{n-1}{n}A_{n-2}' + \frac{8}{n^2\pi} 
= 
\frac{8}{\pi} \left( \frac{1}{n^2} +  \frac{n!}{(n!!)^2} \sum_{{j}=0}^{(n-3)/2}\frac{(2{j}-1)!!}{(2{j})!!}\frac{1}{2{j}+1}  \right).
\]
Now, it is easy to verify that 
\begin{equation}\label{e.easy}
\frac{1}{n^2} = \frac{n!}{(n!!)^2}
\cdot
\frac{(n-2)!!}{(n-1)!!}
\cdot 
\frac{1}{n} 
\end{equation}
and if we substitute this into the previous equation then we get 
\[
\frac{n-1}{n}A_{n-2}' + \frac{8}{n^2\pi} 
= 
\frac{8}{\pi} \left( \frac{n!}{(n!!)^2}
\cdot
\frac{(n-2)!!}{(n-1)!!}
\cdot 
\frac{1}{n} +  \frac{n!}{(n!!)^2} \sum_{{j}=0}^{(n-3)/2}\frac{(2{j}-1)!!}{(2{j})!!}\frac{1}{2{j}+1}  \right),
\]
and we add that first term on the right into the sum (as the ${j}=(n-1)/2$ term) to give us
\[
\frac{n-1}{n}A_{n-2}' + \frac{8}{n^2\pi} 
= 
\frac{8}{\pi} \left( \frac{n!}{(n!!)^2} \cdot \sum_{{j}=0}^{(n-1)/2}\frac{(2{j}-1)!!}{(2{j})!!}\frac{1}{2{j}+1}\right).
\]
Since the right-hand side of the above equation
is the right-hand side of equation 
(\ref{e.bothsidesApprime}), we can re-write the above equation as 
\[
\frac{n-1}{n}A_{n-2}' + \frac{8}{n^2\pi} 
= 
A_n',
\]
thus establishing the validity of  equation (\ref{e.ApApp}), as desired. 

Hence, since both $A_n$ and $A_n'$ satisfy the same recursion from equation (\ref{e.ApApp}), and since they also start at the same value of $A_1= A_1' = 8/\pi$, then $A_n = A_n'$ for all odd values of $n$. This gives us
the desired equality in equation (\ref{e.secondAp}) in the statement of our theorem. 

Next, we will show that the numbers
\begin{equation*}
{n!}  \sum_{{j}=0}^{(n-1)/2} \frac{(2{j}-1)!!}{(2{j})!!} \frac{1}{2{j}+1}
\end{equation*}
from equation (\ref{e.A296726}) really are the same as  every other entry in \rm \seqnum{A296726}, which is the list of coefficients for the exponential generating function for $\arcsin x/(1-x)$.  To show this, we begin with the series for $1/(1-x)$ which is
\begin{equation}\label{e.11minusx}
\frac{1}{1-x} = 1 + x + x^2 + x^3 + x^4 + x^5 + \cdots,
\end{equation}
and for $\arcsin x$ which is
\[
\arcsin x = x + \frac{1}{3!}x^3 + \frac{9}{5!}x^5 + \frac{225}{7!} x^7 + \frac{11025}{9!}x^9 + \cdots,
\]
thanks to \seqnum{A177145}. And furthermore, thanks to 
\seqnum{A001818}, we can re-write those numerators in the above equation as follows:
\[
\arcsin x = \frac{((-1)!!)^2}{1!}x + \frac{(1!!)^2}{3!}x^3 + \frac{(3!!)^2}{5!}x^5 + \frac{(5!!)^2}{7!} x^7 + \frac{(7!!)^2}{9!}x^9 + \cdots.
\]
Hence, since the generating function for 
$\arcsin x/(1-x)$ will be the convolution of the generating functions for 
$\arcsin x$ and $1/(1-x)$, then the
$n$th term in the {\em exponential} generating function for $\arcsin x/(1-x)$,
for $n$ odd,
will be 
\[
n! \left( \frac{((-1)!!)^2}{1!} + \frac{(1!!)^2}{3!} + \frac{(3!!)^2}{5!} + \cdots + \frac{((n-2)!!)^2}{n!}\right),
\]
which we can write as 
\[
n! \sum_{{j}=0}^{(n-1)/2} \frac{((2{j}-1)!!)^2}{(2{j}+1)!}.
\]
Now, since $(2{j}+1)! = (2{j}-1)!!(2{j})!!(2{j}+1)$, then the above expression becomes
\[
n! \sum_{{j}=0}^{(n-1)/2} \frac{(2{j}-1)!!}{(2{j})!!} \frac{1}{2{j}+1},
\]
as seen in equation  (\ref{e.A296726}).

{\bf We conclude with $n$ even}. We  now define $A_n''$ to be the right-hand side of equation (\ref{e.secondApeven}), so  that 
\begin{equation}\label{e.bothsidesApprime.even}
 A_n'' =  \frac{16}{\pi}\cdot \frac{n!}{(n!!)^2} \cdot \sum_{{j}=0}^{(n-2)/2} \frac{(2{j})!!}{(2{j}+1)!!} \frac{1}{2{j}+2}.
\end{equation}
We know from Theorem \ref{t.recursion} that 
for $n$ even, we have 
\[
A_n = \frac{n-1}{n}A_{n-2} + 
\frac{16}{n^2\pi}.
\]
We now seek to
prove that 
\begin{equation}\label{e.ApApp.even}
A_n'' = \frac{n-1}{n}A_{n-2}'' + 
\frac{16}{n^2\pi}
\end{equation}
This, along with the fact that $A_2 = A_2'' = 4/\pi$, is all we will need to conclude that $A_n = A_n''$ for all
$n$, thus proving the validity of equation 
(\ref{e.secondApeven}). 

Now, going back to equation  (\ref{e.bothsidesApprime.even}), we replace $n$ with $n-2$ to give us 
\[
A_{n-2}''
= 
\frac{16}{\pi} \cdot \frac{(n-2)!}{((n-2)!!)^2} 
\cdot
\sum_{{j}=0}^{(n-4)/2}\frac{(2{j})!!}{(2{j}+1)!!}\frac{1}{2{j}+2}.
\]
Next, we multiply both sides by $(n-1)/n$ and then add $16/(n^2\pi)$ to give us 
\[
\frac{n-1}{n}A_{n-2}'' + \frac{16}{n^2\pi}
= 
\frac{16}{\pi} \left( \frac{1}{n^2} + \frac{n-1}{n}\cdot \frac{(n-2)!}{((n-2)!!)^2} 
\cdot
\sum_{{j}=0}^{(n-4)/2}\frac{(2{j})!!}{(2{j}+1)!!}\frac{1}{2{j}+2}\right).
\]
Thanks to equation (\ref{e.easybefore}), the right-hand side simplifies to give us
\[
\frac{n-1}{n}A_{n-2}'' + \frac{16}{n^2\pi}
= 
\frac{16}{\pi} \left( \frac{1}{n^2} + \frac{n!}{(n!!)^2} 
\sum_{{j}=0}^{(n-4)/2}\frac{(2{j})!!}{(2{j}+1)!!}\frac{1}{2{j}+2}\right).
\]
We now use equation (\ref{e.easy}) to replace the $1/n^2$ inside the parenthesis on the right to get
\[
\frac{n-1}{n}A_{n-2}'' + \frac{16}{n^2\pi}
= 
\frac{16}{\pi} \left( 
\frac{n!}{(n!!)^2}\cdot \frac{(n-2)!!}{(n-1)!!} \cdot \frac{1}{n}
+ 
\frac{n!}{(n!!)^2} 
\sum_{{j}=0}^{(n-4)/2}\frac{(2{j})!!}{(2{j}+1)!!}\frac{1}{2{j}+2}\right),
\]
and we add that first term on the right into the sum (as the $j = (n-2)/2$ term) to give us
\[
\frac{n-1}{n}A_{n-2}'' + \frac{16}{n^2\pi}
= 
\frac{16}{\pi} \left(  
\frac{n!}{(n!!)^2} 
\sum_{{j}=0}^{(n-2)/2}\frac{(2{j})!!}{(2{j}+1)!!}\frac{1}{2{j}+2}\right).
\]
Since the right-hand side of the above equation is the right-hand side of equation 
(\ref{e.bothsidesApprime.even}), 
we can re-write the above equation as 
\[
\frac{n-1}{n}A_{n-2}'' + \frac{16}{n^2\pi}
= A_n'',
\]
thus establishing the validity of equation 
(\ref{e.ApApp.even}), as desired. 

Hence, since both $A_n$ and $A_n''$ satisfy the same recursion from equation 
(\ref{e.ApApp.even}), and since they also start at the same value of $A_2 = A_2'' = 4/\pi$, then $A_n = A_n''$ for all even values of $n$. This gives us the desired equality 
(\ref{e.secondApeven})
in the statement of our theorem.

Finally, we will show that the numbers
\begin{equation*}
{n!}  \sum_{{j}=0}^{(n-2)/2} \frac{(2{j})!!}{(2{j}+1)!!} \frac{1}{2{j}+2}
\end{equation*}
from equation (\ref{e.Anotyet}) really are the same as  every other entry in \rm \seqnum{A372324}, which is the list of coefficients for the exponential generating function for 
$\arcsin^2 x/(2(1-x))$.  To show this, 
we begin with the series
\[
\frac{\arcsin x}{\sqrt{1-x^2}} = x 
+ \frac{4}{3!}x^3 + \frac{64}{5!}x^5 + \frac{2304}{7!} x^7 + \frac{147456}{9!}x^9 + \cdots
\]
as seen in \seqnum{A002454}, which also tells us that we can re-write the numerators in 
the above equation as follows:
\[
\frac{\arcsin x}{\sqrt{1-x^2}} = 
\frac{(0!!)^2}{1!}x^1
+ \frac{(2!!)^2}{3!}x^3 + \frac{(4!!)^2}{5!}x^5 + \frac{(6!!)^2}{7!} x^7 + \frac{(8!!)^2}{9!}x^9 + \cdots.
\]
We now integrate both sides to get 
\[
\frac{1}{2}\arcsin^2 x = 
\frac{(0!!)^2}{2!}x^2
+ \frac{(2!!)^2}{4!}x^4 + \frac{(4!!)^2}{6!}x^6 + \frac{(6!!)^2}{8!} x^8 + \frac{(8!!)^2}{10!}x^{10} + \cdots.
\]
Hence, since the generating function for 
$\arcsin^2 x/(2(1-x))$ will be the convolution of equation (\ref{e.11minusx}) and the above equation, then the $n$th term in the 
{\em exponential} generating function 
for $\arcsin^2 x/(2(1-x))$, for $n$ even, will be 
\[
n! \left( \frac{((0)!!)^2}{2!} + \frac{(2!!)^2}{4!} + \frac{(4!!)^2}{6!} + \cdots + \frac{((n-2)!!)^2}{n!}\right),
\]
which we can write as 
\[
n! \sum_{{j}=0}^{(n-2)/2} \frac{((2{j})!!)^2}{(2{j}+2)!}.
\]
Now, since $(2{j}+2)! = (2{j})!!(2{j}+1)!!(2{j}+2)$, then the above expression becomes
\[
n! \sum_{{j}=0}^{(n-2)/2} \frac{(2{j})!!}{(2{j}+1)!!} \frac{1}{2{j}+2},
\]
as seen in equation  (\ref{e.Anotyet}).
\end{proof}

\section{Technical  results}\label{s.technical}

Before we can begin the proofs of Theorems \ref{t.oddpower} and \ref{t.evenpower}, we will need some preliminary results.

\subsection{Three applications of Lagrange's identity}
\begin{lemma}\label{Lagrange}
    Let $N$ and $q$ be positive integers. Then,
\begin{equation*}
    \text{for $N$ even, \ \ \ }
    \sum_{{\ell}=1}^{N/2} \sin \frac{q{\ell}2 \pi}{N} 
    = 
    \begin{cases}
        0,  & \text{for $q$ even;}\\[1.75ex]
		\cot \displaystyle \frac{q \pi}{N},  & \text{for $q$ odd.}
	\end{cases}
\end{equation*}
\end{lemma}
\begin{proof} We call upon Lagrange's trigonometric identity \cite{Balsam}, which states that
\begin{equation}\label{Lagrange1}
\sum_{{\ell}=0}^{m} \sin {\ell} \theta = \frac{\cos \displaystyle \frac{\theta}{2} - \cos\left(m+\displaystyle\frac{1}{2}\right)\theta \vphantom{\bigg|}}{2 \sin \displaystyle\frac{\theta}{2} \vphantom{\bigg|}}.
\end{equation}
Since we are assuming that  $N$ is even, we replace $m$ with $N/2$ and  we replace $\theta$ with $q 2 \pi/N$ in 
equation (\ref{Lagrange1}) to get 
\begin{equation}\label{Lagrange2}
\sum_{{\ell}=0}^{N/2} \sin \frac{q{\ell} 2\pi}{N} = \frac{\displaystyle \cos \frac{q\pi}{N} - \cos\frac{(N+1)q\pi}{N} \vphantom{\bigg|}}
{\displaystyle 2 \sin \frac{q\pi}{N} \vphantom{\bigg|}}.
\end{equation}
Now, $\cos \frac{(N+1)q\pi}{N}$ can be written as 
$\cos \left (q\pi + \frac{q\pi}{N} \right)$, and for $q$ even then $q\pi$ is an even multiple of $\pi$ and so 
$\cos \left (q\pi + \frac{q\pi}{N} \right)$
equals $\cos \frac{q\pi}{N}$.  However, for $q$ odd then 
$q\pi$ is an odd multiple of $\pi$ and so 
$\cos \left (q\pi + \frac{q\pi}{N} \right)$
equals $-\cos \frac{q\pi}{N}$.
When we plug these simplifications into the numerator of equation (\ref{Lagrange2}), we get either 
$0$ or $2 \cos \frac{q\pi}{N}$ in 
the numerator depending on whether $q$ is even or odd, respectively, and 
this gives us our desired formula. 
\end{proof}

\begin{lemma}\label{Lagrange6}
    Let $N$ and $q$ be positive integers. Then,
\begin{equation*}
    \text{for $N$ odd, \ \ \ }
    \sum_{{\ell}=1}^{(N-1)/2} \sin \frac{q{\ell}2 \pi}{N} 
    = 
    \begin{cases}
        \displaystyle -\frac{1}{2} \tan \frac{q \pi}{2 N},  & \text{for $q$ even;}\\[2ex]
		\ \ \, \displaystyle \frac{1}{2} \cot \frac{q \pi}{2 N},  & \text{for $q$ odd.}
	\end{cases}
\end{equation*}
\end{lemma}
\begin{proof} We call upon Lagrange's trigonometric identity 
(\ref{Lagrange1}), 
replacing $m$ with $(N-1)/2$ and 
$\theta$ with $q 2 \pi/N$ 
to give us 
\begin{equation}\label{Lagrange'}
\sum_{{\ell}=0}^{(N-1)/2} \sin \frac{q{\ell} 2\pi}{N} = \frac{\cos \displaystyle \frac{q \pi}{N} - \cos \displaystyle\frac{(N)q\pi}{N}  
\vphantom{\bigg|}}{2 \sin \displaystyle\frac{q\pi}{N} \vphantom{\bigg|}}.
\end{equation}
Now, the last term in the above equation is 
simply $\cos q \pi$, which is either $1$ or $-1$ depending on the parity of $q$. So, 
 the right-hand side of equation (\ref{Lagrange'}) is 
 either $(\cos q\pi/N - 1)/(2 \sin q \pi/N)$ or 
 $(\cos q\pi/N + 1)/(2 \sin q \pi/N)$
depending on whether $q$ is even or odd, respectively. 
And, since $(\cos x - 1)/(2 \sin x) = (-1/2) \tan x/2$ and  $(\cos x + 1)/(2 \sin x) = (1/2) \cot x/2$ by standard trig identities, we have our desired formula. 
\end{proof}

\begin{lemma}\label{Lagrange3}
    Let $N$ and $q$ be positive integers. Then,
\begin{equation*}
    \text{for $N, q$ even, \ \ \ }
    \sum_{{\ell}=1}^{N/2} \sin \frac{q{\ell} \pi}{N} 
    = 
    \begin{cases}
        0,  & \text{for $q \equiv \modd{0} {4}$;}
                    \\[1.75ex]
				\cot \displaystyle \frac{q \pi}{2N},  & \text{for $q \equiv \modd{2} {4}$.}
	\end{cases}
\end{equation*}
\end{lemma}
\begin{proof} We call once more upon Lagrange's trigonometric identity 
(\ref{Lagrange1}).
Since $N$ is even, we will 
 replace $m$ with $N/2$ and  we replace $\theta$ with $q \pi/N$ in 
equation (\ref{Lagrange1}) to get 
\begin{equation}\label{Lagrange3a}
\sum_{{\ell}=0}^{N/2} \sin \frac{q{\ell} \pi}{N} = \frac{\cos \displaystyle \frac{q \pi}{2N} - \cos \displaystyle\frac{(N+1)q\pi}{2N}  
\vphantom{\bigg|}}{2 \sin \displaystyle\frac{q\pi}{2N} \vphantom{\bigg|}}.
\end{equation}
Now,  $\cos\frac{(N+1)q\pi}{2N}$ can be written as 
$\cos \left(\frac{q\pi}{2}+ \frac{q\pi}{2N}\right)$, and for 
$q \equiv \modd{0} {4}$ then $\frac{q\pi}{2}$ is an even multiple of $\pi$  and so  
$\cos \left(\frac{q\pi}{2}+ \frac{q\pi}{2N}\right)$ simplifies to   $\cos \frac{q\pi}{2N}$.
However, for 
$q \equiv \modd{2} {4}$ then 
$\frac{q\pi}{2}$ is an odd multiple of 
$\pi$
is odd
and so 
$\cos \left(\frac{q\pi}{2}+ \frac{q\pi}{2N}\right)$ simplifies to   $-\cos \frac{q\pi}{2N}$
 When we plug these simplifications into
 the right-hand side of equation (\ref{Lagrange3a}), we get either 
$0$ or $2 \cos \frac{q\pi}{2N}$ in the numerator depending on whether $q$ is equivalent to 0 or 2 (mod 4), respectively, and 
this gives us our desired formula. 
\end{proof}

\subsection{A useful limit of cotangents}

\begin{lemma}\label{cotangent}
   For $x$ any real number,
\begin{equation*}
    \lim_{k \to \infty} 
 \left( \frac{1}{k} - 1 \right) \cot \frac{x}{k-1} 
 +
   \left( \frac{1}{k} + 1 \right) \cot \frac{x}{k+1} \ = \ \frac{4}{x}.
\end{equation*}
\end{lemma}
\begin{proof} 
We begin with the Taylor expansion for the cotangent, which gives us
\[
\cot \theta \ = \ \frac{1}{\theta} - \frac{\theta}{3} -\frac{\theta^3}{45} + \cdots 
\ = \ 
\frac{1}{\theta} +  \mathcal{O}(\theta).
\]
If we apply this to our limit, we get
\[
 \left(\frac{1}{k} - 1 \right)
\Bigg(\frac{k-1}{x}  + 
\mathcal{O}\left(\frac{x}{k-1}\right)\Bigg)
\ + \  
\left(\frac{1}{k} + 1 \right)
\Bigg(\frac{k+1}{x}  + 
\mathcal{O}\left(\frac{x}{k+1}\right)\Bigg).
\]
Since $x$ is fixed, we can remove it from inside the $\mathcal{O}$. After expanding the above expression, we get 
\[
\left(\frac{1-k}{k}  \right)
\left(\frac{k-1}{x}\right)  +
\left(\frac{1-k}{k}  \right)\cdot\mathcal{O}\left(\frac{1}{k-1}\right)
\ + \  
\left(\frac{1+k}{k}  \right)
\left(\frac{k+1}{x}\right)  + 
\left(\frac{1+k}{k}  \right)\cdot\mathcal{O}\left(\frac{1}{k+1}\right).
\]
This simplifies nicely to 
\[
\left(\frac{(1-k)(k-1)+ (1+k)(k+1)}{kx }\right) + 
\mathcal{O}\left(\frac{1}{k}\right).
\]
We reduce this to get 
\[
\left(\frac{(1+k)^2 - (1-k)^2}{kx }\right) + 
\mathcal{O}\left(\frac{1}{k}\right) 
\ = \ 
\left(\frac{4k}{kx }\right) + 
\mathcal{O}\left(\frac{1}{k}\right) 
\ = \ \frac{4}{x} + \mathcal{O}\left(\frac{1}{k}\right),
\]
which, as $k \to \infty$, gives us our desired $4/x$. 
\end{proof}

\subsection{Four trig integrals}

\begin{lemma}\label{Bq}
Let $k$ and $q$ be odd numbers. Then, 
if we define
\begin{equation}\label{e.Bq}
C_q = \lim_{k \to \infty}
         \sum_{\ell=1}^{(k-1)/2} \int_{\ell 2\pi/(k+1)}^{\ell 2\pi/(k-1)}
         \Big( \cos qkx - \cos qx\Big) \, dx,
\end{equation}
we have that 
\[ 
C_q = \frac{4}{q^2 \pi}.
\]
\end{lemma}
\begin{proof}
First, we integrate the right-hand side of equation (\ref{e.Bq}) to get 
\begin{equation*}
    C_q = \lim_{k \to \infty}
         \sum_{\ell=1}^{(k-1)/2}
         \left. \frac{1}{kq} \sin qkx - \frac{1}{q} \sin qx \,\,\right|_{x=\ell 2\pi/(k+1)}^{x=\ell 2\pi/(k-1)}.
\end{equation*}
    Taking out the $1/q$ and plugging in the endpoints, we get 
\begin{equation}\label{e.Bq6}
    C_q = \frac{1}{q}\lim_{k \to \infty}
         \sum_{\ell=1}^{(k-1)/2}
         \left( \frac{1}{k} \sin \frac{qk\ell2\pi}{k-1}- \sin \frac{q\ell2\pi}{k-1}\right) 
         - 
         \left( \frac{1}{k} \sin \frac{qk\ell2\pi}{k+1} - \sin \frac{q\ell2\pi}{k+1}\right) 
\end{equation}
Now, if we write
\[
\frac{qk\ell2\pi}{k-1}
= \frac{q(k-1+1)\ell2\pi}{k-1}
= q\ell2\pi + \frac{q\ell2\pi}{k-1}
\]
then we see that 
\begin{equation}\label{e.n-1}
\sin \frac{qk\ell2\pi}{k-1} = \sin \frac{q\ell2\pi}{k-1}.
\end{equation}
Likewise, 
if we write
\[
\frac{qk\ell2\pi}{k+1}
= \frac{q(k+1-1)\ell2\pi}{k-1}
= q\ell2\pi - \frac{q\ell2\pi}{k-1}
\]
then we see that 
\begin{equation}\label{e.n+1}
\sin \frac{qk\ell2\pi}{k+1} = -\sin \frac{q\ell2\pi}{k-1}.
\end{equation}
By substituting equations (\ref{e.n-1}) and (\ref{e.n+1}) into the right-hand side of equation (\ref{e.Bq6}), we have that
\begin{equation*}
     C_q = \frac{1}{q}\lim_{k \to \infty}
         \sum_{\ell=1}^{(k-1)/2}
         \left( \frac{1}{k} - 1 \right) \sin \frac{q\ell2\pi}{k-1}+ \left( \frac{1}{k} + 1 \right) \sin \frac{q\ell2\pi}{k+1}.
\end{equation*}
We now distribute the sum, and change the upper limit of the second summation from $(k-1)/2$ to $(k+1)/2$, which fortunately does
not change the value of the sum, to get 
\begin{equation*}
     C_q = \frac{1}{q}\lim_{k \to \infty}
          \left( \frac{1}{k} - 1 \right)\sum_{\ell=1}^{(k-1)/2}
         \sin \frac{q\ell2\pi}{k-1}
         + \left( \frac{1}{k} + 1 \right) \sum_{\ell=1}^{(k+1)/2}\sin \frac{q\ell2\pi}{k+1}.
\end{equation*}
At this point, since $k$ is odd, then both $k-1$ and $k+1$ are even and so we can apply 
Lemma \ref{Lagrange} (with $q$ odd) to rewrite the above equation as 
\begin{equation*}
     C_q = \frac{1}{q}\lim_{k \to \infty}
          \left( \frac{1}{k} - 1 \right)
          \cot \frac{q \pi}{k-1}
         + \left( \frac{1}{k} + 1 \right) \cot \frac{q \pi}{k+1}.
\end{equation*}
We can now apply Lemma \ref{cotangent} with $x=q\pi$ to the above equation to get that 
\[
     C_q = \frac{1}{q} \frac{4}{q\pi} = \frac{4}{q^2\pi}, 
\]
as desired. 
\end{proof}

\begin{lemma}\label{Cq}
For $k$ odd and $q$ even, 
if we define
\begin{equation}\label{e.Cq}
C_q = \lim_{k \to \infty}
\sum_{\ell=1}^{(k-1)/2} \int_{(\ell-1) \pi/(k-1)}^{\ell \pi/(k+1)}
         f_{q,k}(x) \, dx 
         -
         \int_{\ell \pi/(k+1)}^{\ell \pi/(k-1)}
          f_{q,k}(x) \, dx
\end{equation}
with 
\[
f_{q,k}(x) = \cos qx - \cos qkx,
\]
then we have that 
\[ 
C_q = 
\begin{cases}
    0,  & \text{for $q \equiv \modd{0} {4}$;} \\[1.5ex]
    \displaystyle \frac{16}{q^2 \pi}, & \text{for $q \equiv \modd{2} {4}$.}                    
\end{cases}
\]
\end{lemma}

\begin{proof}
If we let $F_{q,k}(x)$ be the anti-derivative of 
$f_{q,k}(x) = \cos qx - \cos qkx$,
then equation (\ref{e.Cq}) becomes
\begin{equation}\label{e.Cq.secondstage}
C_q = \lim_{k \to \infty}
\sum_{\ell=1}^{(k-1)/2} F_{q,k}(x)
\Bigg|_{(\ell-1)\pi/(k-1)}^{\ell \pi/(k+1)} \, 
 + 
F_{q,k}(x)\Bigg|_{\ell \pi/(k-1)}^{\ell\pi/(k+1)}
\end{equation}
where we replaced $-F_{q,k}$ with $F_{q,k}$ and 
reversed the limits in the second integral. 
We note that almost every term in the above expression for $C_q$ will appear
 twice when we plug in the endpoints and write out the sum,
with the exception of $F_{q,k}(0)$ and $F_{q,k}(\pi/2)$ which will each appear once. 
However, since an easy calculation gives us that 
\begin{equation}\label{e.Fqn.def}
F_{q,k}(x) = \frac{1}{q} \sin qx - \frac{1}{qk} \sin qkx,
\end{equation}
then  $F_{q,k}(0) = 0$ and since
$q$ is even then $F_{q,k}(\pi/2) = 0$ as well.

	So, if we plug in the endpoints, write out the sum, and replace 
the $F_{q,k}(0)$ term with 
$F_{q,k}(\pi/2)$, then equation (\ref{e.Cq.secondstage}) becomes 
\begin{equation*}
C_q =2 \lim_{k \to \infty} \sum_{\ell=1}^{(k-1)/2}  F_{q,k}\left(\frac{\ell\pi}{k+1}\right) - F_{q,k}\left(\frac{\ell\pi}{k-1}\right). 
\end{equation*}
Replacing $F_{q,k}$ with the expression in 
equation (\ref{e.Fqn.def}) and taking out the $1/q$ gives us 
\begin{equation}\label{e.Cqsin}
C_q =\frac{2}{q} \lim_{k \to \infty} \sum_{\ell=1}^{(k-1)/2}  \left( \sin \frac{q\ell\pi}{k+1} - \frac{1}{k} \sin \frac{qk\ell\pi}{k+1} \right) - \left( \sin \frac{q\ell\pi}{k-1} - \frac{1}{k} \sin \frac{qk\ell\pi}{k-1} \right).
\end{equation}
Now, if we write
\begin{equation}\label{e.C.n-1.early}
\frac{qk\ell\pi}{k+1}
= \frac{q(k+1-1)\ell\pi}{k+1}= q\ell\pi - \frac{q\ell\pi}{j+1}
\end{equation}
and if we remember that $q$ is even, then we see that 
\begin{equation}\label{e.C.n-1}
\sin \frac{qk\ell\pi}{k+1} = -\sin \frac{q\ell\pi}{k+1}.
\end{equation}
Likewise, 
if we write
\begin{equation}\label{e.C.n+1.early}
\frac{qk\ell\pi}{k-1} 
= \frac{q(k-1+1)\ell\pi}{k-1} = q\ell\pi + \frac{q\ell\pi}{k-1}
\end{equation}
and again recall that $q$ is even, then we see that 
\begin{equation}\label{e.C.n+1}
\sin \frac{qk\ell\pi}{k-1} = \sin \frac{q \ell \pi}{k-1}.
\end{equation}
By substituting equations (\ref{e.C.n-1}) and (\ref{e.C.n+1}) into the right-hand side of equation (\ref{e.Cqsin}), we have that
\begin{equation*}
     C_q = \frac{2}{q}\lim_{k \to \infty}
         \sum_{\ell=1}^{(k-1)/2}
         \left( 1+\frac{1}{k} \right) \sin \frac{q\ell\pi}{k+1}- \left( 1-\frac{1}{k}  \right) \sin \frac{q\ell\pi}{k-1}.
\end{equation*}
We now distribute the sum, factor through the negative in the second expression, and change the upper limit of the first summation from $(k-1)/2$ to $(k+1)/2$, which fortunately does
not change the value of the sum, to get 
\begin{equation*}
     C_q = \frac{2}{q}\lim_{k \to \infty}
          \left( \frac{1}{k}+1 \right)\sum_{\ell=1}^{(k+1)/2}
         \sin \frac{q\ell\pi}{k+1}
         + \left( \frac{1}{k} -1\right) \sum_{\ell=1}^{(k-1)/2}\sin \frac{q\ell\pi}{k-1}.
\end{equation*}
At this point, since $k$ is odd, then both $k+1$ and $k-1$ are even and so we can apply 
Lemma \ref{Lagrange3} (with $q$ even). If $q/2$ is even, then Lemma \ref{Lagrange3} tells us that both the above sums are zero and so 
$C_q = 0$ in this case. If $q/2$ is odd, we 
apply Lemma \ref{Lagrange3} to tell us that 
\begin{equation*}
     C_q = \frac{2}{q}\lim_{k \to \infty}
          \left( \frac{1}{k} +1\right) 
          \cot \frac{q\pi}{2(k+1)}
         + \left( \frac{1}{k} -1\right) 
         \cot \frac{q\pi}{2(k-1)}
    \qquad \text{for $q/2$ odd.}            
\end{equation*}
We can now apply Lemma \ref{cotangent} with $x=q\pi/2$ to the above equation to get that 
\[
     C_q = \frac{2}{q} \frac{4}{q\pi/2} = \frac{16}{q^2\pi} \qquad \text{for $q/2$ odd,}            
\]
as desired. 
\end{proof}

\begin{lemma}\label{Cq'}
For $k$ even and $q$ odd, 
if we define
\begin{equation}\label{e.Cq'}
C_q = \lim_{k \to \infty}
\sum_{\ell=1}^{k/2} \int_{(\ell-1) 2\pi/(k-1)}^{\ell 2\pi/(k+1)}
         f_{q,k}(x) \, dx 
         -
         \int_{\ell 2\pi/(k+1)}^{\ell 2\pi/(k-1)}
          f_{q,k}(x) \, dx
\end{equation}
such that  
\[
   f_{q,k}(x) = \cos qx - \cos qkx,
\]
then we have
\[
   C_q = \frac{8}{q^2 \pi}.
\]
\end{lemma}

\begin{proof}
If we let $F_{q,k}(x)$ be the anti-derivative of 
$f_{q,k}(x) = \cos qx - \cos qkx$,
then equation (\ref{e.Cq'}) becomes
\begin{equation}\label{e.Cq'.secondstage}
C_q = \lim_{k \to \infty}
\sum_{\ell=1}^{k/2} F_{q,k}(x)
\Bigg|_{(\ell-1)2\pi/(k-1)}^{\ell 2\pi/(k+1)} \, 
 + 
F_{q,k}(x)\Bigg|_{\ell 2\pi/(k-1)}^{\ell2\pi/(k+1)}
\end{equation}
where we replaced $-F_{q,k}$ with $F_{q,k}$ and 
reversed the limits in the second integral. 
We note that almost every term in the above expression for $C_q$ will appear
 twice when we plug in the endpoints and write out the sum,
with the exception of $F_{q,k}(0)$ and 
$F_{q,k}(k \pi/(k-1))$ which will each appear once. 
However, from equation (\ref{e.Fqn.def}) we note that 
 $F_{q,k}(0) = 0$, and as for 
$F_{q,k}(k \pi/(k-1))$, we note that we can use the same technique as in equations (\ref{e.C.n+1.early}) and (\ref{e.C.n+1}), this time with $k$ in place of 
$\ell$, to give us
\[
F_{q,k}(k \pi/(k-1)) = \frac{1}{q} \sin \frac{q k \pi}{k-1} - \frac{1}{qk} \sin \frac{q k \pi}{k-1}.
\]
Since $k$ goes to infinity, we see that 
$F_{q,k}(k \pi/(k-1))$ will vanish. 

	So, returning to 
equation (\ref{e.Cq'.secondstage}), 
if we plug in the endpoints and  write out the sum, 
eliminating the term 
$F_{q,k}(0)$ and doubling the term 
$F_{q,k}(k \pi/(k-1))$ (this is legitimate since both terms will vanish as $k \to \infty$), 
we have
\begin{equation*}
C_q =2 \lim_{k \to \infty} \sum_{\ell=1}^{k/2}  F_{q,k}\left(\frac{\ell 2\pi}{k+1}\right) - F_{q,k}\left(\frac{\ell 2\pi}{k-1}\right). 
\end{equation*}
Replacing $F_{q,k}$ with the expression in 
equation (\ref{e.Fqn.def}) and taking out the $1/q$ gives us 
\begin{equation}\label{e.Cq'sin}
C_q =\frac{2}{q} \lim_{k \to \infty} \sum_{\ell=1}^{k/2}  \left( \sin \frac{q\ell 2\pi}{k+1} - \frac{1}{k} \sin \frac{qk\ell 2\pi}{k+1} \right) - \left( \sin \frac{q\ell 2\pi}{k-1} - \frac{1}{k} \sin \frac{qk\ell 2\pi}{k-1} \right).
\end{equation}
Using the same techniques as seen in equations 
(\ref{e.C.n-1.early}), 
(\ref{e.C.n-1}), 
(\ref{e.C.n+1.early}), and
(\ref{e.C.n+1}), we can write 
\[
\sin \frac{q k \ell 2 \pi}{k+1} = - \sin \frac{q \ell 2 \pi}{k+1} \qquad \text{and} \qquad 
\sin \frac{q k \ell 2 \pi}{k-1} =  \sin \frac{q \ell 2 \pi}{k-1}.  
\]
By substituting these two equations  into the right-hand side of equation (\ref{e.Cq'sin}), we have that
\begin{equation*}
     C_q = \frac{2}{q}\lim_{k \to \infty}
         \sum_{\ell=1}^{k/2}
         \left( 1+\frac{1}{k} \right) \sin \frac{q\ell 2\pi}{k+1}- \left( 1-\frac{1}{k}  \right) \sin \frac{q\ell 2\pi}{k-1}.
\end{equation*}
We now distribute the sum, and we  factor through the negative in the second expression,  to get 
\begin{equation*}
        C_q = \frac{2}{q}\lim_{k \to \infty} \left(
          \left( \frac{1}{k} +1\right) \sum_{\ell=1}^{k/2}
        \sin \frac{q\ell 2\pi}{k+1} + 
        \left( \frac{1}{k}  - 1 \right) \sum_{\ell=1}^{k/2} \sin \frac{q\ell 2\pi}{k-1}\right).
\end{equation*}
We can change the upper limit of the last sum from $k/2$ to $(k-2)/2$, because when $\ell = k/2$ that last term becomes $\sin q \pi k/(k-1)$ which approaches 
$\sin q \pi = 0$ as $k \to \infty$. This gives us 
\begin{equation*}
        C_q = \frac{2}{q}\lim_{k \to \infty}\left(
          \left( \frac{1}{k} +1\right) \sum_{\ell=1}^{k/2}
        \sin \frac{q\ell 2\pi}{k+1} + 
        \left( \frac{1}{k}  - 1 \right) \sum_{\ell=1}^{(k-2)/2} \sin \frac{q\ell 2\pi}{k-1}\right).
\end{equation*}

At this point, since $k$ is even, then both $k+1$ and $k-1$ are odd and so we can apply 
Lemma \ref{Lagrange6} to both of the above sums (the first with $N=k+1$ and the second with $N=k-1$), to give us 
\begin{equation*}
     C_q = \frac{2}{q}\lim_{k \to \infty}\left(
          \left( \frac{1}{k} +1\right) 
          \frac{1}{2}\cot \frac{q\pi}{2(k+1)}
         + \left( \frac{1}{k} -1\right) 
  \frac{1}{2}\cot \frac{q\pi}{2(k-1)}\right). 
\end{equation*}
We can now factor out 1/2 and apply Lemma \ref{cotangent} with $x=q\pi/2$ to the above equation to get that 
\[
     C_q = \frac{1}{q} \left(\frac{4}{q\pi/2}\right) = \frac{8}{q^2\pi},      
\]
as desired. 
\end{proof}

\begin{lemma}\label{Cq*}
For $k$ and $q$ both even, 
if we define
\begin{equation}\label{e.Cq*}
C_q = \lim_{k \to \infty}
\sum_{\ell=1}^{k/2} \int_{(\ell-1) \pi/(k-1)}^{\ell \pi/(k+1)}
         f_{q,k}(x) \, dx 
         -
         \int_{\ell \pi/(k+1)}^{\ell \pi/(k-1)}
          f_{q,k}(x) \, dx
\end{equation}
with 
\[
f_{q,k}(x) = \cos qx - \cos qkx,
\]
then we have that 
\[ 
C_q = 
\begin{cases}
    0,  & \text{for $q \equiv \modd{0} {4}$;} \\[1.5ex]
    \displaystyle \frac{16}{q^2 \pi}, & \text{for $q \equiv \modd{2} {4}$.}                    
\end{cases}
\]
\end{lemma}

\begin{proof}
We proceed as in the proof of Lemma \ref{Cq}. 
If we (again) let $F_{q,k}(x)$ be the anti-derivative of 
$f_{q,k}(x) = \cos qx - \cos qkx$,
then equation (\ref{e.Cq*}) becomes
\begin{equation}\label{e.Cq.secondstage*}
C_q = \lim_{k \to \infty}
\sum_{\ell=1}^{k/2} F_{q,k}(x)
\Bigg|_{(\ell-1)\pi/(k-1)}^{\ell \pi/(k+1)} \, 
 + 
F_{q,k}(x)\Bigg|_{\ell \pi/(k-1)}^{\ell\pi/(k+1)}
\end{equation}
where we replaced $-F_{q,k}$ with $F_{q,k}$ and 
reversed the limits in the second integral. 

We note that almost every term in the above expression for $C_q$ will appear
 twice when we plug in the endpoints and write out the sum,
with the exception of $F_{q,k}(0)$ and $F_{q,k}(k\pi/(2(k-1)))$ which will each appear once. 
From equation (\ref{e.Fqn.def}) we have 
that
\[
F_{q,k}(x) = \frac{1}{q} \sin qx - \frac{1}{qk} \sin qkx,
\]
which tells us that $F_{q,k}(0) =0$. As for 
$F_{q,k}(k\pi/(2(k-1)))$, we have that 
\[
F_{q,k}(k\pi/(2(k-1))) = \frac{1}{q} \sin \frac{qk \pi}{2(k-1)} - \frac{1}{qk} \sin \frac{qk^2 \pi}{2(k-1)}.
\]
As $k \to \infty$, the first term approaches $(1/q) \sin(q \pi/2)$ which is zero since $q$ is even, and the second term approaches zero thanks to the $1/(qk)$ in front. 

	So, if we plug in the endpoints in 
 equation (\ref{e.Cq.secondstage*}), write out the sum, and eliminate 
the terms $F_{q,k}(0)$ 
and 
$F_{q,k}(k\pi/(2(k-1)))$, 
then 
we have 
\begin{equation*}
C_q =2 \lim_{k \to \infty} \sum_{\ell=1}^{k/2}  F_{q,k}\left(\frac{\ell\pi}{k+1}\right) - F_{q,k}\left(\frac{\ell\pi}{k-1}\right). 
\end{equation*}
If we replace $F_{q,k}$ with the expression in 
equation (\ref{e.Fqn.def}) and take out the $1/q$, we get 
\begin{equation*}
C_q =\frac{2}{q} \lim_{k \to \infty} \sum_{\ell=1}^{k/2}  \left( \sin \frac{q\ell\pi}{k+1} - \frac{1}{k} \sin \frac{qk\ell\pi}{k+1} \right) - \left( \sin \frac{q\ell\pi}{k-1} - \frac{1}{k} \sin \frac{qk\ell\pi}{k-1} \right).
\end{equation*}
Since $q$ is even, we can substitute equations 
(\ref{e.C.n-1}) and (\ref{e.C.n+1})
into the above expression to give us
\begin{equation*}
     C_q = \frac{2}{q}\lim_{k \to \infty}
         \sum_{\ell=1}^{k/2}
         \left( 1+\frac{1}{k} \right) \sin \frac{q\ell\pi}{k+1}- \left( 1-\frac{1}{k}  \right) \sin \frac{q\ell\pi}{k-1}.
\end{equation*}
We now distribute the sum, and we factor through the negative in the second expression, to give us 
\begin{equation}\label{e.Cq!!}
     C_q = \frac{2}{q}\lim_{k \to \infty}
        \left(  \left( \frac{1}{k}+1 \right)\sum_{\ell=1}^{k/2}
         \sin \frac{q\ell\pi}{k+1}
         + \left( \frac{1}{k} -1\right) \sum_{\ell=1}^{k/2}\sin \frac{q\ell\pi}{k-1} \right).
\end{equation}

Let us now consider the last term in the last sum of the above equation, namely,
\[
\left( \frac{1}{k} -1\right) \sin \frac{q(k/2)\pi}{k-1}.
\]
As $k \to \infty$, this term approaches
\[
\bigg(  -1\bigg) \sin \frac{q\pi}{2},
\]
and since $q$ is even then this equals zero. Hence, we can safely eliminate the $\ell=k/2$ term from the second sum. Furthermore, using once again that $q$ is even, we can write the above equation (\ref{e.Cq!!}) as 
\begin{equation*}
     C_q = \frac{2}{q}\lim_{k \to \infty}
        \left(  \left( \frac{1}{k}+1 \right)\sum_{\ell=1}^{k/2}
         \sin \frac{(q/2)\ell2\pi}{k+1}
         + \left( \frac{1}{k} -1\right) \sum_{\ell=1}^{(k-2)/2}\sin \frac{(q/2)\ell2\pi}{k-1} \right).
\end{equation*}

We now wish to call upon Lemma \ref{Lagrange6} for the two sums above (the first with $N=k+1$ and the second with $N=k-1$), but we will have to consider the two cases for $q/2$ even and $q/2$ odd. If $q/2$ is even, then
Lemma \ref{Lagrange6} tells us that the above expression is 
\begin{equation*}
\frac{2}{q}\lim_{k \to \infty}
        \Bigg(  \left( \frac{1}{k}+1 \right) 
\left(\frac{-1}{2}\tan \frac{(q/2)\pi}{2(k+1)}\right)
         + \left( \frac{1}{k} -1\right) 
\left(\frac{-1}{2}\tan \frac{(q/2)\pi}{2(k-1)}\right)
\Bigg),
\end{equation*}
and as $k \to \infty$ this is clearly zero. 
If $q/2$ is odd, then from Lemma \ref{Lagrange6} we have 
\begin{equation*}
     C_q = \frac{2}{q}\lim_{k \to \infty}
        \Bigg(  \left( \frac{1}{k}+1 \right) 
\left(\frac{1}{2}\cot \frac{(q/2)\pi}{2(k+1)}\right)
         + \left( \frac{1}{k} -1\right) 
\left(\frac{1}{2}\cot \frac{(q/2)\pi}{2(k-1)}\right)
\Bigg).
\end{equation*}
We now apply Lemma \ref{cotangent} with $x = (q/4) \pi$ to get 
\begin{equation*}
     C_q = \frac{2}{q}\cdot \frac{1}{2} \cdot \frac{4}{(q/4)\pi} = \frac{16}{q^2\pi} \qquad \text{for $q/2$ odd,}
\end{equation*}
as desired. 
\end{proof}

\section{Proof of Theorem \ref{t.oddpower}}\label{s.proof1}

\begin{proof}[Proof of Theorem \ref{t.oddpower}]
We recall that in Theorem \ref{t.oddpower} we are considering the area 
between $\cos^n x$ and 
$\cos^n kx$
where the exponent $n$ is an odd number. We also need to consider the
parity of the coefficient $k$.

{\bf First, we suppose $k$ is odd}. 
As seen in Figure \ref{f.311color} with $n=3$ and $k=11$, there is odd symmetry across the midpoint $x=\pi/2$ and so 
each region ``below" 
$\cos^3 x$ (in color) has an
equivalent area ``above" 
$\cos^3 x$ (in a matching color).

\begin{figure}[h]
\begin{center}
	\includegraphics[width=3.5in]{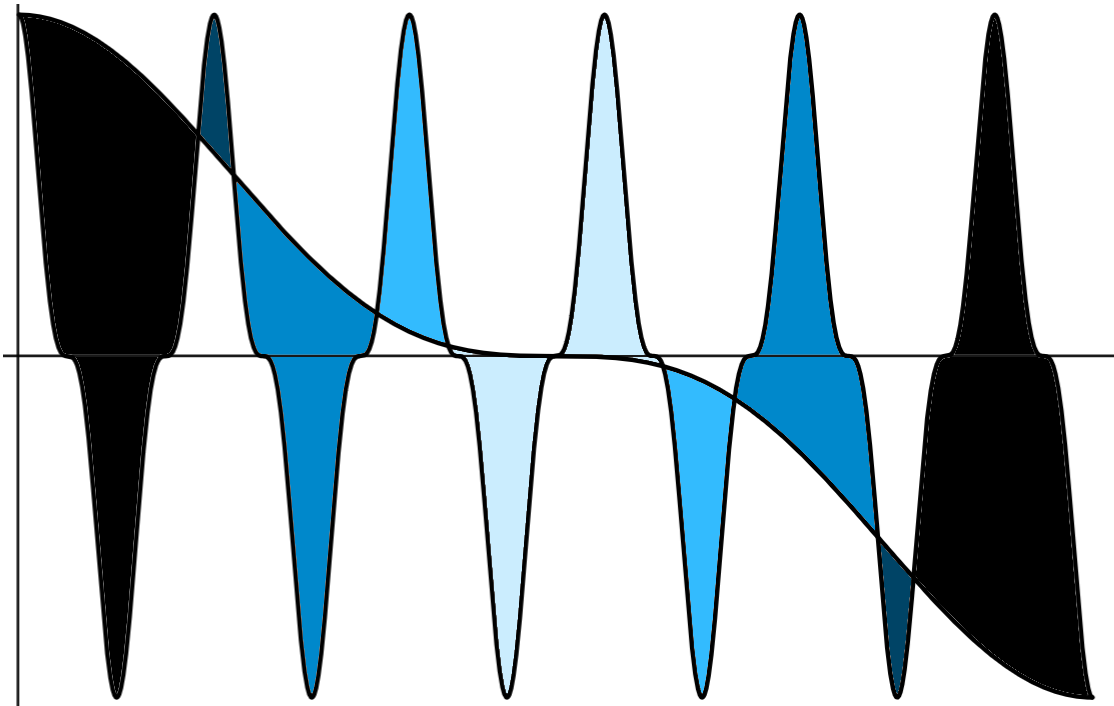}
\end{center}
\caption{Region between $\cos^3 x$ and $\cos^3 11x$.}\label{f.311color}
\end{figure}

In other words, we can just find the areas ``above" $\cos^n x$ on the interval $[0, \pi]$ and then double them. 
To do so, we first need to find the intersection points. 
Since $n$ is odd, then to find the the solutions to 
$\cos^n x = \cos^n kx$ we take the $n$th root of both sides and rewrite it to get
$\cos x - \cos kx = 0$, and we then use a trig identity to write that as 
\[
\sin \frac{(k+1)x}{2}\cdot \sin\frac{(k-1)x}{2} \ = \ 0.
\]
This has solutions $x = \ell\cdot 2\pi /(k+1)$ and 
$x = \ell\cdot 2 \pi /(k-1)$ for $\ell$ any integer, and we note that we can order these as follows:
\begin{multline*}
 0 <  \frac{1\cdot 2\pi}{k+1}
< \frac{1 \cdot 2\pi}{k-1}
<
  \frac{2\cdot 2\pi}{k+1}
< \frac{2 \cdot 2\pi}{k-1}
< \cdots \\[1.3ex]
\cdots  < 
    \frac{\ell\cdot 2\pi}{k+1} 
    <  \frac{\ell\cdot 2\pi}{k-1} 
<
    \frac{(\ell+1)\cdot 2\pi}{k+1} 
    <  \frac{(\ell+1)\cdot  2\pi}{k-1} 
< \cdots
\\[1.3ex]
 \cdots < 
    \frac{(k-1)/2 \cdot 2\pi}{k+1} 
    <  \frac{(k-1)/2 \cdot 2\pi}{k-1} = \pi,
\end{multline*}
and in particular we have that 
\begin{equation}\label{e.inequality}
\frac{\ell\cdot 2\pi}{k-1} 
 <  \frac{(\ell+1) \cdot 2\pi}{k+1}
\qquad \text{so long as $\ell < (k-1)/2$.}
\end{equation}
With these intersection points, we have
the following formula for the total area which takes just the ``upper" regions and doubles them:
\begin{equation}\label{Area.1}
2 \sum_{\ell=1}^{(k-1)/2} \int_{\ell 2\pi/(k+1)}^{\ell 2\pi/(k-1)}
\Big( \cos^n kx - \cos^n x\Big) \, dx.
\end{equation}

We now use the power-reduction formula for cosine to an odd power $n$,
\begin{equation}\label{e.cosoddpowerreduction}
\cos^n \theta = \frac{2}{2^n} \sum_{j=0}^{(n-1)/2} \binom{n}{j} \cos (n-2j)\theta,
\end{equation}
and when we substitute this into equation (\ref{Area.1}), twice, we get the following
expression for the area:
\begin{equation}\label{Area.2}
2 \sum_{\ell=1}^{(k-1)/2} \int_{\ell 2\pi/(k+1)}^{\ell 2\pi/(k-1)}
\frac{2}{2^n} \sum_{j=0}^{(n-1)/2} \binom{n}{j}
\Big( \cos (n-2j)kx - \cos (n-2j)x\Big) \, dx.
\end{equation}
Of course, we want the limit of the expression in  (\ref{Area.2}) as $k$ goes to infinity, 
so when we do this, and re-arrange the sums and integrals and such, we get
\begin{equation}\label{Area.3}
    A_n = \frac{4}{2^n} \sum_{j=0}^{(n-1)/2} 
         \binom{n}{j} \lim_{k \to \infty}
         \sum_{\ell=1}^{(k-1)/2} \int_{\ell 2\pi/(k+1)}^{\ell 2\pi/(k-1)}
         \Big( \cos (n-2j)kx - \cos (n-2j)x\Big) \, dx.
\end{equation}
We now recognize the limit in the
right-hand side of equation (\ref{Area.3}) as
being the same as in Lemma \ref{Bq}, but with the $q$ in that lemma replaced by $n-2j$.
In other words, we now have that 
\begin{align*}
    A_n &= \frac{4}{2^n} \sum_{j=0}^{(n-1)/2} 
         \binom{n}{j} \frac{4}{(n-2j)^2\pi}\\
         &= \frac{8}{\pi} \cdot\frac{1}{2^{n-1} }\cdot
         \sum_{j=0}^{(n-1)/2} 
         \binom{n}{j} \frac{1}{(n-2j)^2}, \label{Apfinal9}
\end{align*}
as desired. 

{\bf Next, we suppose $k$ is even}. 
In contrast to the previous case for $k$ odd, we see in 
Figure \ref{f.310color} with $n=3$ and $k=10$ that when $k$ is even we do 
{\em not} have odd symmetry across the midpoint $x=\pi/2$. Furthermore, the last region on the right is actually cut in half. 

\begin{figure}[h!]
\begin{center}
	\includegraphics[width=3.5in]{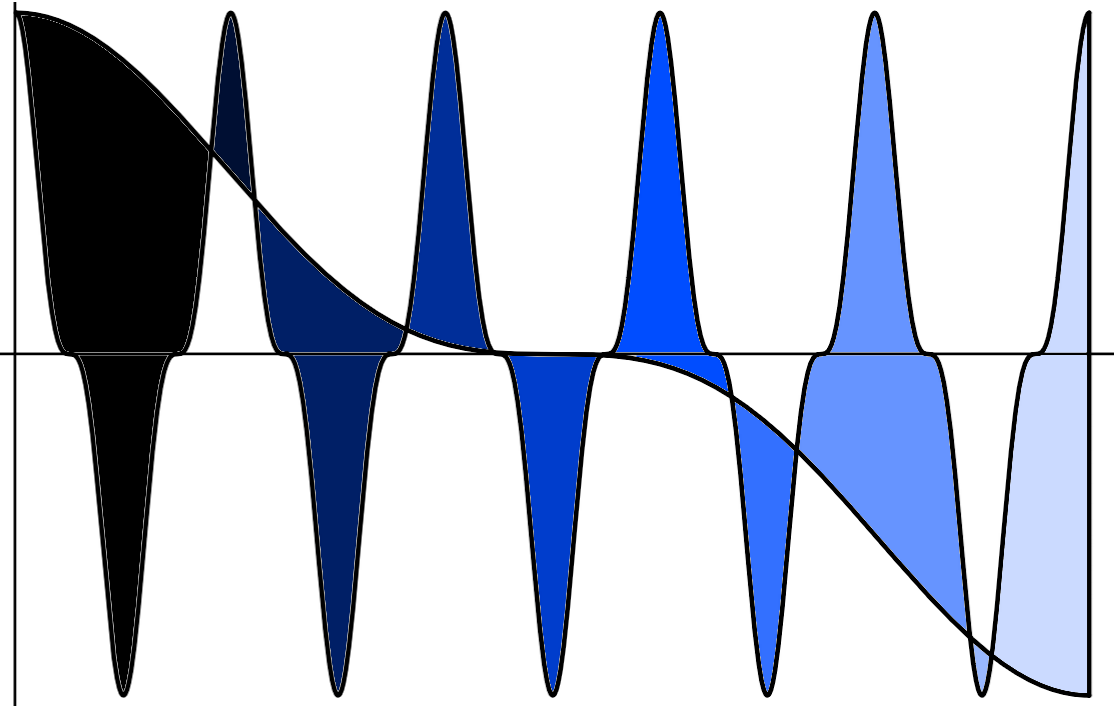}
\end{center}
\caption{Region between $\cos^3 x$ and $\cos^3 10x$.}\label{f.310color}
\end{figure}

Putting all that aside for the moment, we will now find the endpoints for our integrals.
They are mostly the same as the intersection points  but with a slight change towards the end of the list. We have 
\begin{multline*}
 0 < \frac{1\cdot 2\pi}{k+1}
< \frac{1 \cdot 2\pi}{k-1}
<
  \frac{2\cdot 2\pi}{k+1}
< \frac{2 \cdot 2\pi}{k-1}
< \cdots \\[1.3ex]
\cdots <
    \frac{\ell\cdot 2\pi}{k+1} 
    <  \frac{\ell\cdot 2\pi}{k-1} 
<
    \frac{(\ell+1)\cdot 2\pi}{k+1} 
    <  \frac{(\ell+1)\cdot  2\pi}{k-1} 
< 
 \cdots   
    <  \frac{k/2 \cdot 2\pi}{k+1} < \pi,
\end{multline*}
and we note that the last endpoint (namely, 
$\pi$) does not fit the pattern. This is because the last region from Figure 
\ref{f.310color} is cut in half. 

None the less,
we do still have  the inequality from equation (\ref{e.inequality}),
and so 
with all this in mind, we now set up our integrals for the total area. If we set 
$f_{n,k}(x) = \cos^n x - \cos^n kx$, then 
the total area in Figure \ref{f.310color},
moving from left to right, 
is the integral of $f_{n,k}(x)$ over 
$[0, 1\cdot 2\pi/(k+1)]$, followed by the
integral of $-f_{n,k}(x)$
over  
$[1\cdot 2\pi/(k+1), 1\cdot 2\pi/(k-1)]$,
followed by 
the
integral of $f_{n,k}(x)$
over 
$[1\cdot 2\pi/(k-1), 2\cdot 2\pi/(k+1)]$,
and so on. 
The last region should be 
from $(k/2) 2\pi/(k+1)$ to $\pi$, but for convenience we will
consider this as being from 
$(k/2) 2\pi/(k+1)$ to 
$(k/2) 2\pi/(k-1)$ and then subtract
away that ``extra" area from 
$\pi$ to $(k/2) 2\pi/(k-1)$.

In total, then, we have 
\begin{equation*}
\sum_{\ell=1}^{k/2} \left(\int_{(\ell-1)2\pi/(k-1)}^{\ell 2\pi/(k+1)} f_{n,k}(x)\, dx 
+ 
\int_{\ell2\pi/(k+1)}^{\ell 2\pi/(k-1)} -f_{n,k}(x)\, dx\right)
- 
\int_{\pi}^{k\pi/(k-1)} -f_{n,k}(x)\, dx,
\end{equation*}
but of course we actually want the 
{\em limit} of the above expression as $k$ goes to infinity. Since the last term in the above expression is the integral of a bounded function over an interval of width $\pi/(k-1)$, then as $k$ goes to infinity this last term will vanish. So, we have 
\begin{equation}\label{Area.Anafterlimit}
A_n = \lim_{k \to \infty} \sum_{\ell=1}^{k/2} \left(\int_{(\ell-1)2\pi/(k-1)}^{\ell 2\pi/(k+1)} f_{n,k}(x)\, dx 
+ 
\int_{\ell2\pi/(k+1)}^{\ell 2\pi/(k-1)} -f_{n,k}(x)\, dx\right).
\end{equation}

We now use our power-reduction formula from equation (\ref{e.cosoddpowerreduction}) 
to re-write $f_{n,k}(x)$ as
\[
f_{n,k}(x) = \frac{2}{2^n} \sum_{j=0}^{(n-1)/2} \binom{n}{j} 
\Big( \cos (n-2j)x - \cos (n-2j)kx \Big).
\]
When we substitute this into equation (\ref{Area.Anafterlimit}) twice, and distribute the outer sum, we have the following equation for $A_n$:
\begin{multline}\label{Area.Anafterlimit2}
A_n = \frac{2}{2^n} 
\sum_{j=0}^{(n-1)/2} \binom{n}{j} \lim_{k \to \infty} 
\Bigg(\sum_{\ell=1}^{k/2} \int_{(\ell-1)2\pi/(k-1)}^{\ell 2\pi/(k+1)} 
\Big(\cos(n-2j)x - \cos(n-2j)kx \Big) \, dx \\
- 
\sum_{\ell=1}^{k/2} \int_{\ell2\pi/(k+1)}^{\ell 2\pi/(k-1)} 
\Big(\cos(n-2j)x - \cos(n-2j)kx \Big) \, dx\Bigg).
\end{multline}

At this point, we recognize the limit on the 
right-hand side of equation (\ref{Area.Anafterlimit2})
as being the same as in Lemma \ref{Cq'} with $k$ even, but with the 
odd value of $q$ in that lemma replaced by the odd number $n-2j$. In other words, we now have that 
\begin{align*}
    A_n &= \frac{2}{2^n} \sum_{j=0}^{(n-1)/2} 
         \binom{n}{j} \frac{8}{(n-2j)^2\pi}\\
         &= \frac{8}{\pi} \cdot\frac{1}{2^{n-1} }\cdot
         \sum_{j=0}^{(n-1)/2} 
         \binom{n}{j} \frac{1}{(n-2j)^2}, \label{Apfinal9}
\end{align*}
as desired. 
\end{proof}

\section{Proof of Theorem \ref{t.evenpower}}\label{s.proof2}

\begin{proof}[Proof of Theorem \ref{t.evenpower}]
We recall that for this theorem, we are considering the area between 
$\cos^n$ and $\cos^n kx$ where the exponent $n$ is an even number. The proofs are slightly different, depending on the parity of the coefficient $k$.

{\bf First, we suppose $k$ is odd}.
As seen in Figure \ref{f.47color} with $n=4$ and $k=7$, 
there is even symmetry across the midpoint $x=\pi/2$
and so each region on the left of 
$x=\pi/2$ (in color) has an equivalent area on the right of 
$x=\pi/2$ (in a matching color).

\begin{figure}[h]
\begin{center}\includegraphics[width=4.25in]{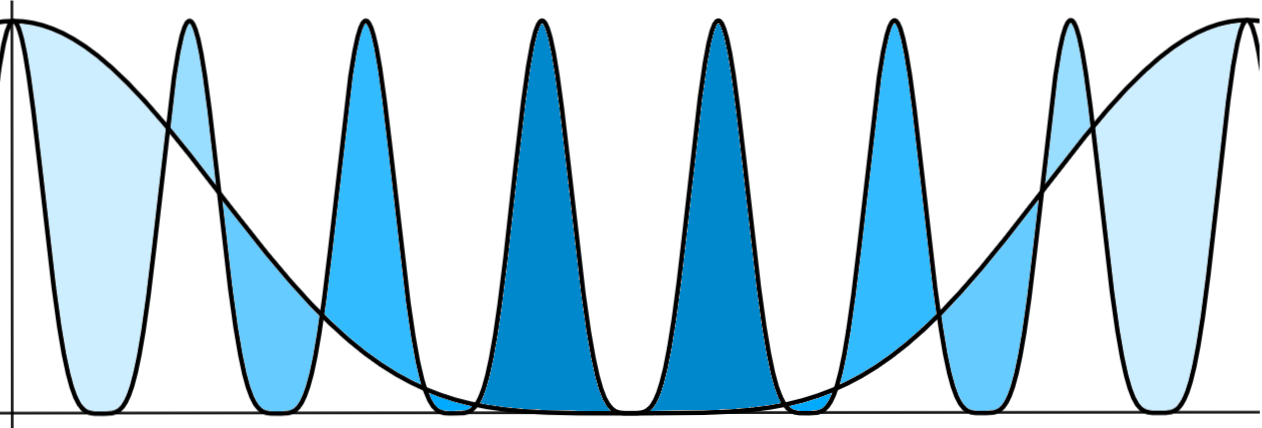}
\end{center}
\caption{Region between $\cos^4 x$ and $\cos^4 7x$.}\label{f.47color}
\end{figure}

In other words, we can just find the 
areas from $0$ to $\pi/2$ and double them. To do so, we first
need to find
the intersection points. 
If we set 
\begin{equation}\label{e.cospeven}
\cos^n x = \cos^n kx
\end{equation}
and take the $n$th root of both sides, then since $n$ is even we will get 
\[
\cos x = \pm \cos kx,
\]
which becomes two equations, 
\[
\cos x - \cos kx = 0 \qquad \text{and} \qquad \cos x + \cos kx = 0.
\]
Using two familiar trig identities, these become
\[
\sin \frac{(k+1)x}{2}\sin \frac{(k-1)x}{2} = 0 \qquad \text{and} \qquad 
	\cos \frac{(k+1)x}{2}\cos \frac{(k-1)x}{2} = 0 
\]
The first equation has solutions
 $x=0$, and also $x = 2\pi /(k+1)$ and 
$x = 2\pi /(k-1)$, and also 
$x = 4\pi /(k+1)$ and 
$x = 4\pi /(k-1)$, and so on. 
The second equation has solutions 
 $x = \pi /(k+1)$ and 
$x = \pi /(k-1)$,
and also 
 $x = 3\pi /(k+1)$ and 
$x = 3\pi /(k-1)$, and so on. 
Hence, the complete list of solutions to equation (\ref{e.cospeven}) in the 
interval $[0, \pi/2]$, written out
in order, is
\begin{multline*}
 0   <   \frac{\pi}{k+1}
     <   \frac{\pi}{k-1}
     <   \frac{2\pi}{k+1}
     <   \frac{2\pi}{k-1}
     <   \frac{3\pi}{k+1}
     <   \frac{3\pi}{k-1}
     <   \cdots \\
   			\cdots 
   				  <  \frac{(k-1)/2\cdot \pi}{k+1}
   				  <   \frac{(k-1)/2\cdot \pi}{k-1}  
   				  =    \frac{\pi}{2}. 
\end{multline*}
With these intersection points,
we have the following formula for
the total area (for $k$ any fixed odd number) which takes just the 
regions on the left of $x=\pi/2$ and doubles them:
\begin{equation}\label{Area.p4}
2 \sum_{\ell=1}^{(k-1)/2} \int_{(\ell-1)\pi/(k-1)}^{\ell \pi/(k+1)} f_{n,k}(x)\, dx + 
\int_{\ell\pi/(k+1)}^{\ell \pi/(k-1)} -f_{n,k}(x)\, dx,
\end{equation}
where $f_{n,k}(x) = \cos^n x - \cos^n kx$.

We now use the power-reduction formula for cosine to an even power $n$,
\[
\cos^n \theta = \frac{1}{2^n} \binom{n}{n/2} + \frac{2}{2^n} \sum_{j=0}^{(n/2)-1} \binom{n}{j} \cos (n-2j)\theta,
\]
to give us that 
\[
f_{n,k}(x) = \frac{2}{2^n} \sum_{j=0}^{(n/2)-1} \binom{n}{j} 
\Big( \cos (n-2j)x - \cos (n-2j)kx \Big).
\]
When we substitute this into equation (\ref{Area.p4}) twice, and distribute the outer sum, we get the following
expression for the area:
\begin{multline}\label{Area.2p42}
2 \sum_{\ell=1}^{(k-1)/2} \int_{(\ell-1) \pi/(k-1)}^{\ell \pi/(k+1)}
\frac{2}{2^n} \sum_{j=0}^{(n/2)-1} \binom{n}{j}
\Big( \cos (n-2j)x - \cos (n-2j)kx\Big) \, dx
\\ 
- 2 \sum_{\ell=1}^{(k-1)/2} \int_{\ell \pi/(k+1)}^{\ell \pi/(k-1)}
\frac{2}{2^n} \sum_{j=0}^{(n/2)-1} \binom{n}{j}
\Big( \cos (n-2j)x - \cos (n-2j)kx\Big) \, dx. 
\end{multline}
Of course, we want the limit of the expression in  (\ref{Area.2p42}) as $k$ goes to infinity, 
so when we do this, and re-arrange the sums and integrals and such, we get
\begin{multline}\label{Area.3p4}
    A_n = \frac{4}{2^n} \sum_{j=0}^{(n/2)-1}
         \binom{n}{j} \lim_{k \to \infty}
         \Bigg( \sum_{\ell=1}^{(k-1)/2} \int_{(\ell-1)\pi/(k-1)}^{\ell \pi/(k+1)}
         \Big( \cos (n-2j)x - \cos (n-2j)kx\Big) \, dx \\
         - 
         \sum_{\ell=1}^{(k-1)/2} \int_{\ell \pi/(k+1)}^{\ell \pi/(k-1)}
         \Big( \cos (n-2j)x - \cos (n-2j)kx\Big) \, dx \Bigg)
\end{multline}
We now recognize the limit 
in the right-hand side of equation (\ref{Area.3p4}) as being the same as in Lemma \ref{Cq}. In other words, we now have that
\begin{equation}    \label{e.Appenultimate}
A_n = \frac{4}{2^n} \sum_{j=0}^{(n/2)-1}
         \binom{n}{j}C_{n-2j}
\end{equation}
where $C_{n-2j}$ from Lemma \ref{Cq} is
defined as 
\[ 
C_{n-2j} = 
 \begin{cases}
        0,  & \text{for $n-2j \equiv \modd{0} {4}$;}\\[1.5ex]
	    \displaystyle \frac{16}{(n-2j)^2 \pi}, & \text{for $n-2j \equiv \modd{2} {4}$.}
    \end{cases}
\]

We now consider the case when 
$n \equiv \modd{2} {4}$. In this case, if we write out the terms in equation (\ref{e.Appenultimate}) and use our definition of $C_{n-2j}$ from above, we have only the terms with $j$ even (as that is when $n - 2j \equiv \modd{2} {4}$), giving us
\[
A_n = \frac{4}{2^n} 
\left(         \binom{n}{0}\frac{16}{(n)^2 \pi}
+ \binom{n}{2}\frac{16}{(n-4)^2\pi} 
+ \binom{n}{4}\frac{16}{(n-8)^2\pi} + \cdots 
+ \binom{n}{(n/2)-1}\frac{16}{(2)^2\pi}
\right)
\]
We now factor out $16/(2^2 \pi)$ from each term, giving us 
\[
A_n = \frac{4}{2^n} \frac{16}{2^2 \pi}
\left(         \binom{n}{0}\frac{1}{(n/2)^2}
+ \binom{n}{2}\frac{1}{(n/2-2)^2} 
+ \binom{n}{4}\frac{1}{(n/2-4)^2} + \cdots 
+ \binom{n}{(n/2)-1} \frac{1}{(1)^2}
\right)
\]
We re-index the above sum, and simplify the coefficients on the left, to get
\[
A_n = \frac{16}{\pi} \cdot \frac{1}{2^{n}} \cdot \sum_{j=0}^{(n-2)/4} \binom{n}{2j} \frac{1}{( n/2 - 2j)^2},
\]
as desired (for $n \equiv \modd{2} {4}$). 

Finally, for $n \equiv \modd{0} {4}$, we again 
write out the terms in equation (\ref{e.Appenultimate}) and use our definition of $C_{n-2j}$ from above. This time, the only non-zero contributions come from $j$ odd (as this is when $n - 2j \equiv \modd{2} {4}$), giving us  
\[
A_n = \frac{4}{2^n} 
\left(         \binom{n}{1}\frac{16}{(n-2)^2 \pi}
+ \binom{n}{3}\frac{16}{(n-6)^2\pi} 
 + \cdots 
+ \binom{n}{(n/2)-1}\frac{16}{(2)^2\pi}
\right)
\]
We again factor out $16/(2^2 \pi)$ from each term, giving us 
\[
A_n = \frac{4}{2^n} \frac{16}{2^2 \pi}
\left(     
  \binom{n}{1}\frac{1}{(n/2-1)^2}
+ \binom{n}{3}\frac{1}{(n/2-3)^2} 
 + \cdots 
+ \binom{n}{(n/2)-1}\frac{1}{(1)^2\pi}
\right)
\]
We re-index the above sum, and simplify the coefficients on the left, to get
\[
A_n = \frac{16}{\pi} \cdot \frac{1}{2^{n}} \cdot \sum_{j=0}^{(n-4)/4} \binom{n}{2j+1} \frac{1}{( n/2 - (2j+1))^2},
\]
as desired (for $n \equiv \modd{0} {4}$).

{\bf Next, we suppose $k$ is even}.
As seen in Figure \ref{f.48blue} with $n=4$ and $k=8$, 
we again have even symmetry across the midpoint $x=\pi/2$
and so each region on the left of 
$x=\pi/2$ (in color) has an equivalent area on the right of 
$x=\pi/2$ (in a matching color). The difference from
the earlier part when $k$ was odd is that here, with $k$ even, the ``center" spike needs to be considered as a special case.

\begin{figure}[h]
\begin{center}\includegraphics[width=4.25in]{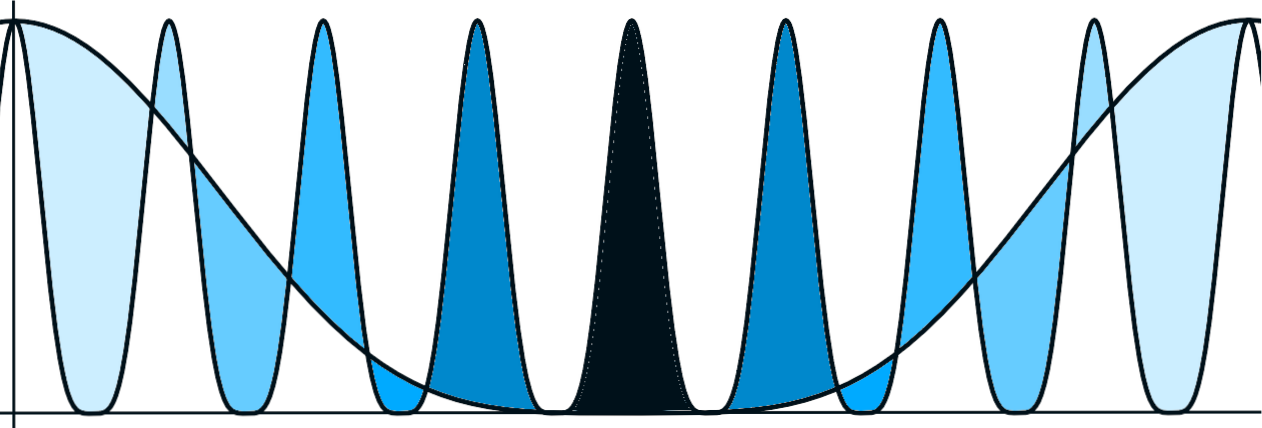}
\end{center}
\caption{Region between $\cos^4 x$ and $\cos^4 8x$.}\label{f.48blue}
\end{figure}

Proceeding as before, the intersection points are mostly the same, except that here we will go just a bit beyond the midpoint of $x=\pi/2$ because that center spike (in black) extends just a bit beyond our interval of $[0,\pi/2]$. As $k$ goes to $\infty$, this extra little bit of area just beyond $x=\pi/2$ will be negligible and so can be ignored. With this in mind, we will integrate over the following sub-intervals:
\[
 0   <   \frac{\pi}{k+1}
     <   \frac{\pi}{k-1}
     <   \frac{2\pi}{k+1}
     <   \frac{2\pi}{k-1}
    <
   			\cdots 
   				  <  \frac{k/2\cdot \pi}{k+1}
   				  <   \frac{k/2\cdot \pi}{k-1}.
\]
With these intersection points, we have the following formula for the total area.
\begin{equation}\label{Area.p4'}
2 \sum_{\ell=1}^{k/2} \int_{(\ell-1)\pi/(k-1)}^{\ell \pi/(k+1)} f_{n,k}(x)\, dx + 
\int_{\ell\pi/(k+1)}^{\ell \pi/(k-1)} -f_{n,k}(x)\, dx,
\end{equation}
where, as before, $f_{n,k}(x) = \cos^n x - \cos^n kx$.

The only difference between equation (\ref{Area.p4'}) above, and equation (\ref{Area.p4}) earlier in the proof (for the case when $k$ odd) is that the upper limit for the sum in  equation (\ref{Area.p4'})  is $k/2$ instead of $(k-1)/2$. So, following the same steps as before, we can jump ahead to what would be our version of equation (\ref{Area.3p4}), which would be 
\begin{multline*}
    A_n = \frac{4}{2^n} \sum_{j=0}^{(n/2)-1}
         \binom{n}{j} \lim_{k \to \infty}
         \Bigg( \sum_{\ell=1}^{k/2} \int_{(\ell-1)\pi/(k-1)}^{\ell \pi/(k+1)}
         \Big( \cos (n-2j)x - \cos (n-2j)kx\Big) \, dx \\
         - 
         \sum_{\ell=1}^{k/2} \int_{\ell \pi/(k+1)}^{\ell \pi/(k-1)}
         \Big( \cos (n-2j)x - \cos (n-2j)kx\Big) \, dx \Bigg).
\end{multline*}
At this point, we recognize the limit in the right-hand side of the above equation as being the same as in Lemma \ref{Cq*}. In other words, we now have that 
\begin{equation*}
A_n = \frac{4}{2^n} \sum_{j=0}^{(n/2)-1}
         \binom{n}{j}C_{n-2j}
\end{equation*}
where $C_{n-2j}$ from Lemma \ref{Cq*} is
defined as 
\[ 
C_{n-2j} = 
 \begin{cases}
        0,  & \text{for $n-2j \equiv \modd{0} {4}$;}\\[1.5ex]
	    \displaystyle \frac{16}{(n-2j)^2 \pi}, & \text{for $n-2j \equiv \modd{2} {4}$.}
    \end{cases}
\]
This is identical to  equation (\ref{e.Appenultimate}) and so we can arrive at the same conclusion here (with $k$ even) as we did before (with $k$ odd). 
\end{proof}

\bigskip
\hrule
\bigskip

\noindent 2020 {\it Mathematics Subject Classification}: Primary 05A10; 
Secondary 05A15, 26A06, 26A42. 

\noindent \emph{Keywords: }  Double factorial, cosine, arcsine, generating function,  binomial coefficient.

\bigskip
\hrule
\bigskip

\noindent (Concerned with sequences
\seqnum{A000129},
\seqnum{A001818},
\seqnum{A002454},
\seqnum{A006882},
\seqnum{A177145},
\seqnum{A296726},
\seqnum{A359311},
and
\seqnum{A372324}.)

\bigskip
\hrule
\bigskip


\begin{thebibliography}{10}
 
\bibitem{Balsam}
Jonathan Balsam, 
Proof Without Words: Lagrange’s Trigonometric Identity (Part II), {\em College Math J.} {\bf 54} (2023), 235.


\bibitem{2191} John Chapman, Yvonne Cheng, and 
Greg Dresden,  Problem \#2191. \textit{Math. Mag.} 
{\bf 97} (2024), 223.


\bibitem{Gould}
Henry Gould and 
Jocelyn Quaintance, 
Double Fun with Double Factorials,
{\em Math. Mag.} {\bf 85} (2012), 177--192.


\bibitem{oeis}
 N.~J.~A.~Sloane et al., The On-Line Encyclopedia of Integer Sequences, 2024. Available at \url{https://oeis.org}.


\end{thebibliography}
\end{document}